\documentclass[11pt]{article} % Thesis form

\usepackage{pdfsync}
\usepackage[centertags]{amsmath}
\usepackage{amsfonts}
\usepackage{amssymb}
\usepackage{amsthm}
\usepackage{newlfont}
\usepackage[utf8]{inputenc}
\usepackage{mathrsfs}
\usepackage{setspace}
\usepackage{fancyhdr}
\usepackage{epsfig}
\usepackage{graphicx}
\usepackage{accents}

\usepackage{makeidx}
\makeindex

\usepackage[footnotesize,bf]{caption}

\usepackage{verbatim} % for \begin-end{comment} block comments

\DeclareMathSymbol{\widehatsym}{\mathord}{largesymbols}{"62}

\newtheorem{defn}{Definition}[section]

\newtheorem{example}[defn]{Example}

\newtheorem{lem}[defn]{Lemma}

\newcommand{\pf}{\noindent{\bf Proof }\mbox{   }}

\newcommand{\EV}{\mathbb{E}}    % expected value of a random variable
\newcommand{\Z}{\mathbb{Z}}     % integers
\newcommand{\mc}[1]{\mathcal{#1}}

\newcommand{\ZK}{\mathbb{Z}^2_K} % 2-dim torus of size K
\newcommand{\bd}{\partial}	    % boundary of a set

\newcommand{\inv}{^{-1}}

\newcommand{\ie}{\emph{i.e.}, }
\newcommand{\eg}{\emph{e.g.}, }
\newcommand{\st}{such that }

\newcommand{\Dnsc}{D(0,n+s)^c}
\newcommand{\Dnc}{D(0,n)^c}
\newcommand{\Dn}{D(0,n)}
\newcommand{\Dns}{D(0,n+s)}

\newcommand{\hp}{\hat{\pi}_K}
\newcommand{\hDnc}{\hp(\Dnc_K)}
\newcommand{\hDnsc}{\hp(\Dnsc_K)}
\newcommand{\hDn}{\hp(\Dn)}
\newcommand{\hDns}{\hp(\Dns)}
\newcommand{\hdDns}{\hp(\bd D(0,n)_s)}
\newcommand{\hdDnsc}{\hp((\bd D(0,n)_s)^c_K)}

\newcommand{\hDR}{\hp(D(0,R))}

\newcommand{\hDr}{\hp(D(0,r))}

\newcommand{\hTDnc}{T_{\hDnc}}

\newcommand{\hTDn}{T_{\hp(D(0,n))}}
\newcommand{\hTdDns}{T_{\hp(\bd D(0,n)_s)}}
\newcommand{\hTDr}{T_{\hp(D(0,r))}}

\newcommand{\hTDRc}{T_{\hp(D(0,R)^c_K)}}

\newcommand{\hGDnc}{\hat{G}_{\hp(\Dnc_K)}}

\newcommand{\hGDn}{\hat{G}_{\hp(D(0,n))}}

\def\midformat{
\setlength{\itemsep}{0pt} \setlength{\parindent}{0mm}
\setlength{\parskip}{0.12in} \setlength{\textheight}{180mm} % 20110726 to make spaces between paragraphs larger
\setlength{\textwidth}{150mm} \setlength{\evensidemargin}{0in}
\setlength{\oddsidemargin}{0in} \setlength{\topmargin}{0in}
\setlength{\hoffset}{1.0cm} \setlength{\voffset}{0.0cm}
\setlength{\headheight}{15pt} \setlength{\headsep}{.5in}
\setlength{\headwidth}{150mm} } \midformat

\newtheoremstyle{theorem}% name
 {}% Space above
 {}%            Space below
 {\itshape}%    Body font
 {}%            Indent amount (empty = no indent, \parindent = para indent)
 {\ttfamily}%   Thm head font
 {.}%           Punctuation after thm head
 {.5em}%        Space after thm head: " " = normal interword space;
 {}%            Thm head spec (can be left empty, meaning `normal')

\newtheoremstyle{plaintext}% name
 {}% Space above
 {}% Space below
 {\upshape}% Body font
 {}% Indent amount (empty = no indent, \parindent = para indent)
 {\ttfamily}% Thm head font
 {.}% Punctuation after thm head
 {.5em}% Space after thm head: " " = normal interword space;
 {}% Thm head spec (can be left empty, meaning `normal')

\theoremstyle{theorem}

\theoremstyle{plaintext}

\numberwithin{equation}{section}
\begin{document}
\makeatother
%%%%%%%%%%%%%%%%%%%%%%%%%
%   Introduction    %
%%%%%%%%%%%%%%%%%%%%%%%%%
\pagenumbering{roman} \thispagestyle{empty}
\title{On Escaping, Entering, and Visiting Discs \\of Projections of Planar Symmetric Random Walks \\on the Lattice Torus}
%\title{On Escaping, Entering, and Visiting Balls \\of Projections of Symmetric Random Walks \\on the Periodic Lattice}
\author{Michael Carlisle\footnote{michael.carlisle@baruch.cuny.edu}\\
Baruch College, CUNY}
\maketitle

\begin{abstract}
We examine escape and entrance times, Green's functions, local times, and hitting distributions of discs and annuli 
of a symmetric random walk on $\Z^2$ projected onto the periodic lattice $\Z^2_K$. This extends a framework for the simple planar random walk in \cite{DPRZ2006} to the large class of planar random walks in \cite{BRFreq}.
The approach uses comparisons between $\Z^2$ and $\Z^2_K$ hitting times and distributions on annuli, and uses only random walk methods. 
\end{abstract}

%-------------

\oddsidemargin 0.0in \textwidth 6.0in \textheight 8.5in

\singlespacing
%\setcounter{tocdepth}{0} % only print ch in toc, no sec or subsec #ing 
%\setcounter{tocdepth}{2} % print ch.sec.subsec in toc
%\tableofcontents 

%\setcounter{tocdepth}{2} % only print ch.sec in toc, no subsec #ing 

%  \cleardoublepage 
 % \addcontentsline{toc}{chapter}{\listfigurename} 
  %\listoffigures 
%\newpage

%  \cleardoublepage 
%\addcontentsline{toc}{chapter}{\indexname}
%\printindex   % why is this breaking ALL formatting??!!

%\doublespacing  %ZZZ UNCOMMENT THIS FOR FINAL VERSION
\pagenumbering{arabic}
\pagestyle{fancy}
\renewcommand{\headrulewidth}{0.0pt}
%\headwidth
\lhead{} \chead{} \rhead{\thepage}
\fancyhead[RO]{\thepage} \fancyfoot{}

\section{Introduction} \label{ch:Intro}

There is a wealth of literature on random walks on the planar lattice $\Z^2$: 
 Aldous (\cite{AF}), Dembo, Peres, Rosen, \& Zeitouni (\cite{DPRZ2001}, \cite{DPRZ2004}, \cite{DPRZ2006}), Lawler (\cite{LawInt}, \cite{LawCov}, \cite{LawPol}), and Rosen (\cite{RosenET}) all discuss problems of the simple random walk on $\Z^2$; in \cite{BRFreq}, Rosen \& Bass extend certain results to a class of infinite-range symmetric random walks. This paper builds on these works, to examine the timing structure of entrances to and escapes from discs in $\Z^2$, projected onto the square lattice torus $\Z^2_K$.

Consider a random walk $S_t = S_0 + \sum_{j=0}^t X_j$, \index{s@$S_t$} for $X = \{X_j\}_{j \in \mathbb{N} \cup \{0\}}$ \index{x@$X_j$} with the following properties: $S$ is symmetric, $X_1$ has finite covariance matrix equal to a scalar times the identity, \ie $\Gamma := cov(X_1) = cI$, $c>0$, and $X$ is strongly aperiodic.\footnote{\cite{BRFreq} requires, for walks in $\Z^2$, that the covariance matrix of $X_1$ be equal to $\frac{1}{2} I$, but this is a convenience for three technical points (on pages 9, 12, and 42), relating only to rotations. 

It is worthy (if not elementary) to note that the simple random walk on $\Z^d$'s $X_1$ covariance matrix is cov$(X_1) = \frac{1}{d}I$. If $K$ is odd, this walk projects to a strongly aperiodic simple random walk on $\Z^d_K$.} Set $\pi_{\Gamma} := 2\pi \sqrt{\det \Gamma}$\footnote{For symmetric simple random walk on $\Z^d$, $\Gamma = \frac{1}{d} I$, so $\det \Gamma = d^{-d}$; for $d=2$, $\pi_{\Gamma} = \pi$.}. \index{pigamma@$\pi_{\Gamma}$} 

$X_1$ has, for some $\beta>0$ and $M := 4 + 2\beta$, \index{M@$M$}
\begin{equation} \label{eqn:moments}
\EV|X_1|^{M} = \sum_{x \in \Z^2} |x|^{M} p_1(x) < \infty,
\end{equation}
where, as usual in the literature,
\begin{align*}
p_1(x,y) = p_1(y-x) = P^x(X_1 = y)
\end{align*}
is the one-step transition probability, \index{p@$p_1$} and $P^x$ is the probability measure for walks starting at $x$.
The random walk methods used in this paper require $M > 4$ to make escape results on the lattice torus $\Z^2_K$ look as they do on the plane in $\Z^2$.

We will switch between the planar and toral lattice representations of the random walk and corresponding stopping times, hitting distributions, etc. 
Define the projections, for $x = (x_1, ..., x_d) \in \Z^d$, by \index{pi@$\pi_K$} \index{pihat@$\hp$} 
\[\begin{array}{ll}
\pi_K: \Z^d \to [-K/2,K/2)^d \cap \Z^d, \\
\pi_K(x) = \left((x_1 + \lfloor \frac{K}{2} \rfloor) (\text{mod } K) - \lfloor \frac{K}{2} \rfloor, ..., (x_d + \lfloor \frac{K}{2} \rfloor) (\text{mod } K) - \lfloor \frac{K}{2} \rfloor\right); \\
 \hp: \Z^d \to \Z^d_K, 
\,\,\,\, \hp(x) = (\pi_K x) + (K \Z)^d.
\end{array}\]
(For example, if $d=2$, $x = (-12,6)$ and $K = 11$, then $\pi_{11}(\Z^2) = \{-5,\ldots,5\}^2$, $\pi_{11}(x) = (-1,-5)$, and $\hat \pi_{11}(x) = (-1,-5) + (11\Z)^2$.)

\begin{figure}[!ht]
  \centering
    \includegraphics[width=6in]{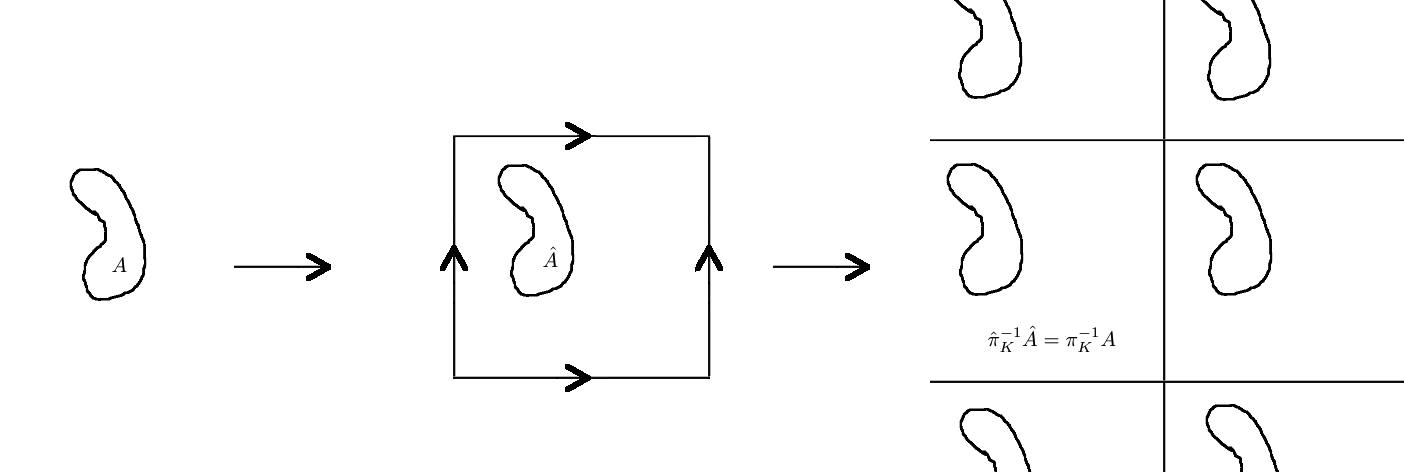}
  \parbox{4in}{
  \caption{$A \to \hat{A} \to \hp\inv \hat{A} = \pi_K \inv A$}
  \label{fig:plane_torus_pullback}}
\end{figure}

We call the set of lattice points $\pi_K(\Z^d) = [-K/2,K/2)^d \cap \Z^d$ the {\bf primary copy} \index{primary copy} in $\Z^d$, and for $x \in \pi_K(\Z^d)$, $\hat{x} := \hp x$ is its corresponding element in $\Z^d_K$. Any $z \in \pi_K \inv x$, $z \neq \pi_K x$, is called a {\bf copy} \index{copy} of $x$. Likewise, for a set $A \subset \Z^d$, $\hat{A} := \hp A$ is the periodic projection of $A$, and the set of all copies of $A$ is \index{piinv@$\hp\inv \hat{A}$}
\[\pi_K \inv\pi_K A = \hp\inv \hat{A} := \{z \in \Z^d: z = x + (j_1 K, ..., j_d K), \,\,  j_i \in \Z, i=1,...,d, x \in A\}.\]
Figure \ref{fig:plane_torus_pullback} displays, for $d=2$, the projection of a planar set $A$ onto the torus as $\hat{A}$, and its pullback onto $\pi_K \inv A$. (If $A \subset \pi_K \Z^d$, then of course, $A = \pi_K A$.)

For a given $\hat{x} \in \Z^d_K$, we define $x := \pi_K \hp \inv \hat x$ to be the primary copy of that element.

While $X_j$ is the $j$th step of the walk and $S_j$ its position at time $j$, we use $\hat{S}_j$ \index{X@$\hat{S}_j$} to denote the position of the projected walk at time $j$. The distance between two points $x,y \in \Z^d$ will be the Euclidean distance $|x-y|$; on the torus, the distance between two points $\hat{x},\hat{y} \in \Z_K^d$ \index{xhat@$\hat{x}$} will be the minimum Euclidean distance $|\hat{x}-\hat{y}|$. To limit the issues regarding this distance, we will restrict any discs on $\Z^d_K$ to have radius $n < K/4$ (sometimes written as a diameter constraint: $2n < K/2$).

To bound our functions, we need a precise notion of bounding distance on the lattice torus $\Z^d_K$. As in \cite{DPRZ2006}, a function $f(x)$ is said to be $O(x)$ if $f(x)/x$ is bounded, uniformly in all implicit geometry-related quantities (such as $K$). That is, $f(x) = O(x)$ if there exists a universal constant $C$ (not depending on $K$) such that $|f(x)| \leq Cx$. Thus $x = O(x)$ but $Kx$ is \emph{not} $O(x)$. A similar convention applies to $o(x)$.

Next, we will define a few terms describing the distance of a random walk step, relative to a reference disc of radius $n$ and an $s$-sized annulus around the disc.
A {\bf small} jump \index{jump}%!small}
 refers to a step that is short enough to possibly (but not necessarily) stay inside a disc of radius $n$ (\ie $|X_1| < 2n$).
A {\bf baby} jump %\index{jump!baby} 
refers to a small jump that is too short to hop over an $s$-annulus from inside a disc (\ie $|X_1| < s$).
A {\bf medium} jump %\index{jump!medium} 
refers to a step that is sufficiently large to hop out of a disc and past an $s$-annulus, but with magnitude strictly less than $K$, and cannot land near a projected copy of its launching point (\ie $s < |X_1| < K-2n$).
A {\bf large} jump %\index{jump!large} 
is a step which, in the projection, would be considered ``wrapping around'' in one step (\ie $|X_1| > K-2n$).
A {\bf targeted} jump %\index{jump!targeted} 
is a large jump which lands directly in a copy of the disc or annulus just launched from. These terms will aid in dealing with differences between regular and projected hitting and escape times.%\footnote{We have distinguished between three types of jumps on the torus that in the planar-only case (as in \eg \cite{BRFreq}) are referred to only as large jumps.}

The paper is structured as follows. 
In Section \ref{ch:Escape}, we prove results about probabilities of exiting a disc in $\Z^2$ and $\Z^2_K$. 
Section \ref{ch:Entry} contains results involving entering a disc. 
In Section \ref{ch:Annulus}, we use the general framework from \cite{CarThree} for analyzing moving between three sets that partition a sample space, and discuss the application of these ideas to hitting an annulus just outside a disc, and gambler's ruin estimates in that case. 

%END Ch1

%\include{Ch2} % results related to exit times
%\section{Escape} % Disc escape functions 
\section{Disc Escape} \label{ch:Escape} % Disc escape functions 

In this section we develop the notions of hitting time and Green's function on $\Z^2$ and $\Z^2_K$, and find relationships between the two with respect to the timing of the random walk's escape from a disc.

\subsection{Disc escape time}
%\subsection{Hitting and escape time} 
\label{sec:HittingEscapeTime}

The \emph{hitting time} of a random walk to a set $A$ is defined as the stopping time \index{time@$T_A$} $T_{A} = \inf\{t \geq 0: S_t \in A\}$. Likewise, the \emph{escape time} of the walk from $A$ is the stopping time $T_{A^c}$. For a recurrent, strongly aperiodic, irreducible random walk on $\Z^2$, $T_{A^c} < \infty$ a.s. We denote $T_{\hat{A}}$ to be the hitting time of $\hat{A} \subset \Z_K^2$. We will examine several relationships between planar and toral hitting times.

An immediate observation on hitting times (\eg from \cite{Spitzer}) is that, the larger the set to hit, the quicker it will be hit. If $A \subset B$, then obviously $T_B \leq T_A$. It is clear, then, that $\hp\inv \hat{A}$, as an infinite number of copies of $A \subset \Z^2$, has a quicker hitting time than just one copy of $A$. In fact, we have 
\begin{equation} \label{eq:PTtimeLaw}
T_{\pi_K \inv A} = T_{\hp\inv \hat{A}} = T_{\hat{A}}.
\end{equation}

\begin{figure}[!ht]
  \centering
    \includegraphics[width=6in]{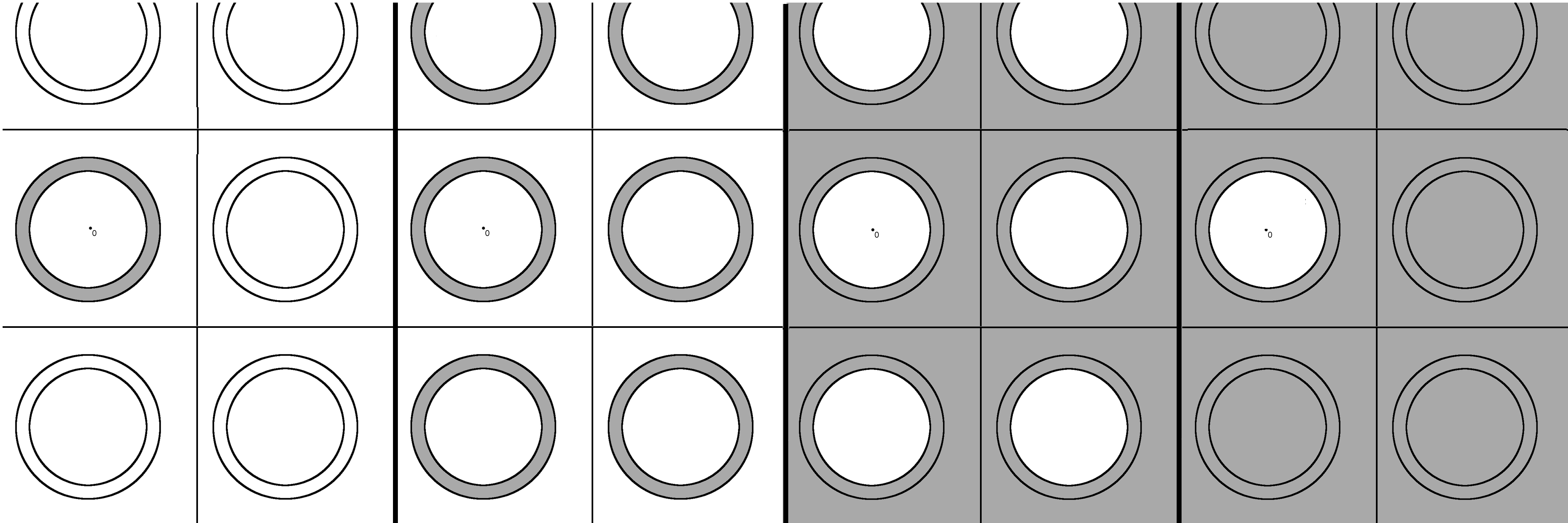}
    \parbox{6in}{
    \caption{Comparison of planar sets listed in \eqref{eq:AnnDiscSubs}, on the plane. Labeled sets are shaded.}
    \label{fig:discs_annuli}}
\end{figure}

Let $n, s$ be such that $n+s < K/4$, and $D(0,n) = \pi_K \Dn$ the primary copy of $\Dn \subset \Z^2$. Define the primary copy's portion of the complement of $\Dn$ to be $\Dnc_K := \Dnc \cap \pi_K \Z^2$. \index{DncK@$\Dnc_K$}
\eqref{eq:AnnDiscSubs} and Figure \ref{fig:discs_annuli} describe the nestedness of sets from the planar annulus $\bd D(0,n)_s$ up to the planar disc complement $\Dnc$: 
\begin{align}
\bd D(0,n)_s & \subset & \pi_K \inv (\bd D(0,n)_s) = \hp\inv \hdDns \notag\\
 & \subset & \hp\inv \hDnc = \pi_K \inv (\Dnc_K) & \subset D(0,n)^c. \label{eq:AnnDiscSubs}
\end{align}

By \eqref{eq:PTtimeLaw}, \eqref{eq:AnnDiscSubs} yields, starting at any $x \in D(0,n)$, the \emph{disc escape time inequalities}
\begin{align} 
T_{\bd D(0,n)_s} & \geq T_{\pi_K \inv \bd D(0,n)_s} = T_{\hp\inv \hdDns} \notag \\
 & \geq T_{\hp\inv \hDnc} = T_{\pi_K \inv (\Dnc_K)} \geq T_{D(0,n)^c} \geq 1. \label{eq:HittingTimeComp1}
\end{align}

%Without loss of generality, w
We shall take planar starting points from the primary copy ($x = \pi_K x$). The probabilities of these inequalities being strict (\eg $P^x(T_{D(0,n)^c} < \hTDnc)$) and the means of the stopping times will be of interest to us. We start with estimating the mean of the planar escape time from $D(0,n)$ (which improves on \cite[Prop. 6.2.6]{LawSRW}), and then use this probability to estimate the toral escape time from $\hDn$.

\begin{lem} \label{lem:EscapeDisc}
Let $S_t = S_0 + \sum_{j=1}^t X_j$ be a random walk in $\Z^2$ with $E|X_1|^2 < \infty$, and covariance matrix $\Gamma$ such that $tr(\Gamma) = \gamma^2 > 0$. Then,
uniformly for $x \in D(0,n)$, and for sufficiently large $n$, \index{expected@$\EV^x(T_{D(0,n)^c})$}
\begin{equation} \label{eqn:EscapeDiscExp}
\frac{n^2 - |x|^2}{\gamma^2} \leq \EV^x(T_{D(0,n)^c}) \leq \frac{n^2 - |x|^2}{\gamma^2} + 2n + 1.
\end{equation}
\end{lem} % also, now see \cite[Prop. 6.2.6, p. 152]{LawSRW}

\pf By \cite[Exercise 1.4]{LawSRW}, the process $M_t := |S_t|^2 - \gamma^2 t$ is a martingale. 

For any given $k$, $k \land T_{D(0,n)^c}$ is a bounded stopping time, and $T_{D(0,n)^c} < \infty$ a.s., so by the monotone convergence theorem, 
\begin{align}
\lim_{k \to \infty} \EV^x(k \land T_{D(0,n)^c}) = \EV^x(T_{D(0,n)^c}). \label{eq:MCTescape}
\end{align}
Hence, by the optional stopping theorem, uniformly for $x \in D(0,n)$,
\begin{equation} \label{eqn:OST_Mt}
 \EV^x(M_{k \land T_{D(0,n)^c}}) = \EV^x(M_0) = |x|^2.
\end{equation}

Decompose $|S_{k \land T_{D(0,n)^c}}|^2$ along the time $k$: 
\begin{align}
|S_{k \land T_{D(0,n)^c}}|^2 & = 1_{\{k \geq T_{D(0,n)^c}\}}|S_{T_{D(0,n)^c}}|^2
  + 1_{\{k < T_{D(0,n)^c}\}}  |S_{k}|^2. \label{eq:decomposeS}
\end{align}
Its expectation, then, is 
\begin{align}
E(|S_{k \land T_{D(0,n)^c}}|^2) & = E\left( 1_{\{k \geq T_{D(0,n)^c}\}}|S_{T_{D(0,n)^c}}|^2   \right)
 + E \left(   1_{\{k < T_{D(0,n)^c}\}}  |S_{k}|^2   \right). \label{eq:decomposeES}
\end{align}
Then by the MCT again, since $T_{D(0,n)^c}<\infty$ a.s.,
\begin{align}
\lim_{k \to \infty} E\left( 1_{\{k \geq T_{D(0,n)^c}\}}|S_{T_{D(0,n)^c}}|^2  \right) =  E\left( |S_{T_{D(0,n)^c}}|^2  \right).  \label{eq:decompES2}
\end{align}
For the second term, note that $1_{\{k < T_{D(0,n)^c}\}}  |S_{k}|^2  \leq n^{2}$, and also 
$1_{\{k < T_{D(0,n)^c}\}}  |S_{k}|^2 \to 0$ a.s. since, again, $T_{D(0,n)^c}<\infty$ a.s.
Thus by the dominated convergence theorem, 
\begin{align}
\lim_{k \to \infty} E\left(   1_{\{k < T_{D(0,n)^c}\}}  |S_{k}|^2   \right) = 0. \label{eq:decompES3}
\end{align}

Combining \eqref{eq:MCTescape}-\eqref{eq:decompES3} yields the expected time 
\begin{align}
|x|^2 & = \EV^x(M_{T_{D(0,n)^c}}) = \EV^x(|S_{T_{D(0,n)^c}}|^2) - \gamma^2 \EV^x(T_{D(0,n)^c}) \notag\\
 \implies \EV^x(T_{D(0,n)^c}) & = \frac{\EV^x(|S_{T_{D(0,n)^c}}|^2) - |x|^2}{\gamma^2}. \label{eq:combined} 
\end{align}
We can bound $|S_{k \land T_{D(0,n)^c}}|^2$ by decomposing along its escape jump: 
if $X_j = (X_j^{(1)}, X_j^{(2)})$ is the $j$th step, then for any $j \leq T_{D(0,n)^c}$, 
\begin{align}
|S_{j}|^2 & = |S_{j-1}|^2 +2   S_{j-1} \cdot X_j + |X_j|^2\leq n^2 + 2n(|X_j^{(1)}| + |X_j^{(2)}|) + |X_j|^2. \label{eqn:St_bounds}
\end{align}
It is clear that $tr(\Gamma) = \gamma^2 = \EV(|X_j|^2)$, and, since $X_j^{(i)} \in \Z$, $i=1,2$, then
\[|X_j^{(1)}| + |X_j^{(2)}| \leq |X_j^{(1)}|^2 + |X_j^{(2)}|^2 = |X_j|^2.\]
Therefore, \eqref{eqn:St_bounds} becomes, substituting $k \land T_{D(0,n)^c} = j \leq T_{D(0,n)^c}$, 
\[ |S_{k \land T_{D(0,n)^c}}|^2 \leq n^2 + (2n+1)|X_{k \land T_{D(0,n)^c}}|^2. \]
By taking expectations,
\begin{align*}
\EV^x(|S_{k \land T_{D(0,n)^c}}|^2) \leq n^2 + (2n+1) \EV^x(|X_{k \land T_{D(0,n)^c}}|^2) \leq n^2 + (2n+1) \gamma^2 < \infty. 
\end{align*}
Since, at $T_{D(0,n)^c}$, we have escaped the disc, we have a lower bound as well. By \eqref{eq:decomposeES}-\eqref{eq:decompES3}, 
\begin{align}
n^2 \leq \EV^x(|S_{T_{D(0,n)^c}}|^2) \leq n^2 + (2n+1)\gamma^2. \label{eq:Sbounds}
\end{align}
Combining \eqref{eq:Sbounds} with \eqref{eq:combined} yields \eqref{eqn:EscapeDiscExp}. \qed

For $\Gamma = cI$, $\gamma^2 = 2c$ and so \eqref{eqn:EscapeDiscExp} becomes\footnote{For simple random walk on $\Z^2$, $c = 1/2$, which yields \cite[(2.3)]{DPRZ2006}.} 
\begin{equation} \label{eqn:EscapeDiscExpCalc}
\frac{n^2 - |x|^2}{2c} \leq \EV^x(T_{D(0,n)^c}) \leq \frac{n^2 - |x|^2}{2c} + 2n + 1.
\end{equation}

We define the \emph{Green's function} \index{Green@$G_{t^*}$} for two points $x, y$, as the expected number of visits to $y$, starting from $x$, up to the fixed time $t^*$:
\begin{equation} \label{eq:GreenDefn}
G_{t^*}(x,y) := \EV^x\bigg[\sum_{j=0}^{t^*} 1_{\{S_j=y\}}\bigg] = \sum_{j=0}^{\infty} P^x(S_j = y; j < t^*).
\end{equation}
Spitzer, in \cite{Spitzer}, similarly defines the \emph{truncated Green's function}, \index{Green@$G_A$} for $x, y \in A$ of a random walk from $x$ to $y$ before exiting $A$ as the total expected number of visits to $y$, starting from $x$:
\begin{equation} \label{eq:GreenDefn}
G_A(x,y) := \EV^x\bigg[\sum_{j=0}^{\infty} 1_{\{S_j=y; j<T_{A^c}\}}\bigg] = \sum_{j=0}^{\infty} P^x(S_j = y; j < T_{A^c})
\end{equation}
and 0 if $x$ or $y \not \in A$. (Since the walk is recurrent and aperiodic, there is no ``all-time'' Green's function to count the total number of visits to $x$ from $j=0$ to $\infty$.) An elementary result for any random walk (found, for example, in \cite{Spitzer}, or \cite[Sect. 1.5]{LawInt}) is that, for $x, y \in A \subset B$, there are more possible visits inside $B$ than inside $A$:
\begin{equation} \label{eq:GreenComp}
 G_A(x,y) \leq G_B(x,y).
\end{equation}
Also of interest is the expected hitting time identity 
\begin{equation} \label{eq:GreenHit}
\EV^x(T_{A^c}) = \sum_{z \in A} G_A(x,z).
\end{equation}
Starting at a point $x \in A^c$, the \emph{hitting distribution} \index{hitting@$H_A$} of $A$ is defined as
\[H_A(x,y) := P^x(S_{T_A} = y).\]
The \emph{last exit decomposition} \index{last exit decomposition} of a hitting distribution is based on the Green's function: for $A$ a proper subset of $\Z^2$, $x \in A^c$, $y \in A$, 
\begin{equation} \label{eq:LastExit}
H_{A}(x,y) = \sum_{z \in A^c} G_{A^c}(x,z) p_1(z,y).
\end{equation}
If $y \in A \subset B$, then for $x \in B^c \subset A^c$, we have by \eqref{eq:GreenComp} the monotonicity result
\begin{equation} \label{eq:HittingComp}
\begin{array}{lll}
H_{A}(x,y) & = & \sum_{z \in A^c} G_{A^c}(x,z) p_1(z,y)\\
 & \geq & \sum_{z \in B^c} G_{B^c}(x,z) p_1(z,y) = H_{B}(x,y)
\end{array}
\end{equation}
and the subset hitting time relations (assuming a recurrent random walk) 
\begin{align} 
P^x(T_A = T_B) & = \sum_{z \in A} H_B(x,z); \notag\\
P^x(T_A \neq T_B) & = P^x(T_A > T_B) = \sum_{z \in B \setminus A} H_B(x,z) \label{eq:SubsetHittingComp}
\end{align}
which we will revisit in Section \ref{ch:Annulus}.

By Markov's inequality, large jumps are rare: if $C_M = \EV(|X_1|^M) < \infty$, then since $2n < K/2$, 
\begin{equation} \label{eqn:LargeJumpEstimate}
P(|X_1|>K-2n) \leq \frac{C_M}{(K-2n)^M} < \frac{2^M C_M}{K^M} = O(K^{-M}).
\end{equation}
Recall that, when given a toral element $\hat{x} \in \ZK$, we define $x$ to be the (planar) primary copy of that element; $x := \pi_K \hp \inv \hat x$. A toral step $\hat{x} \to \hat{y}$ must take into account large jumps that, on the plane, would land on a copy of $y$ (\ie in $\hp\inv \hat{y}$). All of these positions, together, are a small addition to the planar jump probability. By \eqref{eqn:LargeJumpEstimate} we have, for $\hat{x}, \hat{y} \in \Z^2_K$, the targeted jump estimate \index{phat@$\hat{p}_1$}
\begin{align} \label{eqn:TargetedJumpEstimate}
\hat{p}_1(\hat x, \hat y) := P^{\hat{x}}(\hat S_1 = \hat{y}) & = P^{x}(S_1 = y) + P^{x}\left(|X_1|>K-2n ; S_1 \in \hp\inv\hat{y} \setminus \{y\}\right) \notag\\
 & \leq p_1(x, y) + O(K^{-M}).
\end{align}

By \eqref{eq:LastExit}, \eqref{eqn:LargeJumpEstimate}, and then \eqref{eqn:EscapeDiscExp} and \eqref{eq:GreenHit}, for some $c < \infty$ and any $x\in D(0,n)$, 
\begin{align}
P^x(\hTDnc>T_{D(0,n)^c}) & = \sum_{z \in \left(\hp\inv\hDn \,\setminus\, D(0,n)\right)} \sum_{y \in D(0,n)} G_{D(0,n)}(x,y) p_1(y,z) \notag\\
 & \leq c K^{-M} \sum_{y \in D(0,n)} G_{D(0,n)}(x,y) = O(K^{-M} n^2). \label{eqn:ProbTorusExit}
\end{align}

We now find that the mean of the disc escape time on the torus is larger than on the plane, but only by a small factor (induced by the rarity of targeted jumps). \index{expectedT@$\EV^{\hat x}(\hTDnc)$}

\begin{lem} \label{lem:EscapeBounds}
For $n < K/4$, $x \in D(0,n)$, and $n$ and $K$ sufficiently large,  
\begin{equation} \label{eqn:Escape}
\EV^{\hat x}[\hTDnc] \leq \EV^x[T_{D(0,n)^c}] + O(K^{-M} n^2) \max_{y\in D(0,n)}\EV^y[T_{D(0,n)^c}].
\end{equation}
\end{lem}

\pf 
To bound the disc escape time above, consider a ``worst case'' scenario (making the $\hDn$-escape time as long as possible) where every large jump targets the same point inside the disc. 

Let $y^*$ be a point on $D(0,n)$ such that $\EV^{y^*}(T_{D(0,n)^c}) = \max_{y \in D(0,n)} \EV^{y}(T_{D(0,n)^c})$. 
Define the times $\sigma_i$ and $\tau$, and index variable $N$, by 
\begin{align}
\sigma_0 = T_{D(0,n)^c}; \,\, \sigma_{i+1} & = \inf\{j > \sigma_i: S^*_{j-1} + X_j \in D(0,n)^c\}, \,\, i \geq 0 \label{eq:EscapeTimes}\\
\tau & = \inf\left\{ j > 0: |X_j| \leq K-2n, S^*_{j-1} + X_j \in D(0,n)^c \right\} \notag\\
N & = j \iff \sigma_j = \tau \notag
\end{align}
where $\sigma_0$ is the original walk $S$'s planar disc escape time, and the modified walk $S^*$ is defined as the walk whose large jumps (of size $> K-2n$) target $y^*$, until the walk escapes $D(0,n)$ via a nonlarge jump:
\begin{align}
S^*_t := \left\{ \begin{array}{ll}
x, & t=0 \\
x + \sum_{k=1}^t X_k, & 0 < t < \sigma_0 \\
x + \sum_{k=1}^{\tau} X_k, & t \geq \tau \text{ on } \{\tau = \sigma_0\} \\
y^*, & t = \sigma_i, \,\, |X_{\sigma_i}| > K-2n, \,\, 0 \leq i < N \\
y^* + \sum_{k=\sigma_i+1}^t X_k, & \sigma_i < t < \sigma_{i+1}, \,\, |X_{\sigma_i}| > K-2n, \,\, 0 \leq i < N \\
y^* + \sum_{k=\sigma_{N-1}+1}^{\tau} X_k, & t \geq \tau \text{ on } \{ \tau = \sigma_N, \,\, N > 0\}.
\end{array}\right. \label{eq:WorstCaseEscapeWalk}
\end{align}

\begin{figure}[!ht]
  \centering
    \includegraphics[width=4in]{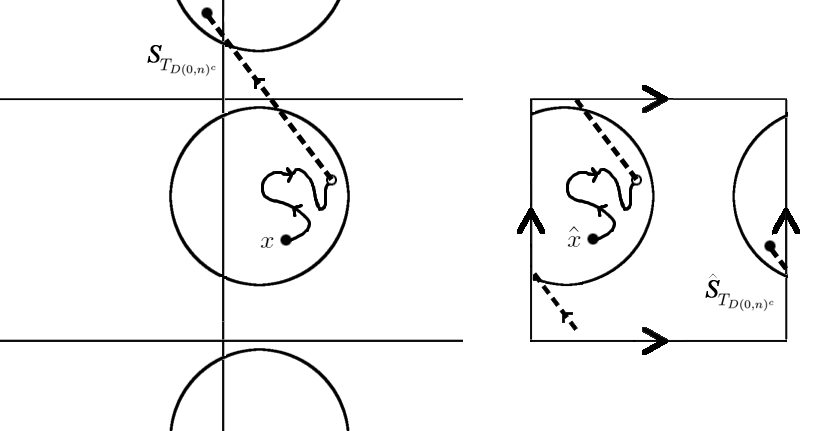}
    \parbox{4in}{
    \caption{An example of a path in $\{\hTDnc > T_{D(0,n)^c}\}$, where a targeted jump of planar distance $\geq O(K)$ keeps the walk in $\hDn$.}
    \label{fig:target_wrap_plane_torus}}
\end{figure}

$\sigma_i$, $i \geq 0$, are the successive would-be escape times from $D(0,n)$, if $y^*$-targeting was not ``enabled''. $\tau$ is the smallest $\sigma_i$ such that escape from $D(0,n)$ actually occurs, and $N$ is the number of large jumps before this escape occurs. Note that, considering times on the original walk $S$, 
\begin{align*}
\{ N = 0 \} & = \left\{ \tau = T_{\hDnc} = T_{\Dnc} \, , \, |X_{T_{\Dnc}}| \leq K-2n \right\} \\
\{ N > 0 \} & = \left\{ \tau \geq T_{\hDnc} > T_{\Dnc} \, , \, |X_{T_{\Dnc}}| > K-2n \right\}, 
\end{align*}
and, conditioned on $\{ N > 0 \}$, $N$ is a geometric random variable with success parameter $p = P^{y^*}\left( T_{D(0,n)^c} = T_{\hp(D(0,n)^c_K)} \, , \, |X_{T_{\Dnc}}| \leq K-2n \right) = 1 + O(K^{-M} n^2)$ by \eqref{eqn:ProbTorusExit} and \eqref{eqn:LargeJumpEstimate} 
(where a ``failure'' is a targeted jump back to $y^* \in \Dn$). Thus, $\EV^{\hat{x}}[\hTDnc] \leq \EV^{\hat{x}}[\tau]$, since $\tau$ is the escape time of $\hDn$, with targeting back to $y^*$. Conditioning on $\{N>0\}$, and by \eqref{eqn:ProbTorusExit} and \eqref{eqn:LargeJumpEstimate} and the strong Markov property on $\sigma_0$, we have 
\begin{align}
\EV^{x}[\tau] & = \EV^{x}[\tau | N=0]P^x(N=0) + \EV^{x}[\tau | N>0] P^x(N>0) \notag \\
 & = \EV^{x}[\sigma_0 | N=0]P^x(N=0) + \EV^{x}[\sigma_0 + \sigma_N - \sigma_0 | N>0] P^x(N>0) \notag \\
 & \leq \EV^{x}[T_{D(0,n)^c}] + \EV^{x}[\sigma_N - \sigma_0 | N>0] P^x(N>0) \label{eq:EscapeExpValUp1} \\
 & \leq \EV^{x}[T_{D(0,n)^c}] + O(K^{-M} n^2) \EV^{y^*}[\sigma_N - \sigma_0 | N>0]. \notag
\end{align}
On $\{N > 0\}$, the time of the $j$th excursion from $y^*$ until attempted disc escape is $\nu_{j} := \sigma_{j} - \sigma_{j-1}$, for $1 \leq j \leq N$, are IID with mean $\EV[\nu_{j}] = \EV^{y^*}[T_{\Dnc}]$.
Since $P(N < \infty) = 1$, by Wald's identity we have 
\begin{align*}
\EV^{y^*}[\sigma_N - \sigma_0 | N>0] & = \EV^{y^*}\left[ \sum_{j=1}^N \nu_j \bigg| N > 0\right] = \EV^{y^*}[N | N>0] \, \EV^{y^*}[\nu_1] \\
 & = \frac{1}{p} \, \EV^{y^*}[T_{\Dnc}] = (1 + O(K^{-M} n^2)) \EV^{y^*}[T_{\Dnc}].
\end{align*}
Therefore, \eqref{eq:EscapeExpValUp1} becomes 
\begin{align*}
\EV^{x}[\tau] & \leq \EV^{x}[T_{D(0,n)^c}] + O(K^{-M} n^2) \EV^{y^*}[\sigma_N - \sigma_0 | N>0] \\
 & \leq \EV^{x}[T_{D(0,n)^c}] + O(K^{-M} n^2) \EV^{y^*}[T_{\Dnc}]. \,\, \qed
\end{align*}

\begin{comment} % BEGIN COMMENTED-OUT EXAMPLE
\begin{example} \label{eg:WorstCaseTiming}
Let $A=D(0,\sqrt{2}) = \{0, +e_1, -e_1, +e_2, -e_2\} \subset \Z^2$, where $e_i$ is the $i$th unit vector in $\Z^2$, and $K$ odd and fixed. Let $X$ be the symmetric random walk on $\Z^2$ starting at $X_0=0$ defined by the probabilities
\[p_1(K^j e_i) = P^0(X_1= K^j e_i) = \frac{1}{4}e^{-\lambda}\frac{\lambda^j}{j!}, j=0,1,2,\ldots; \,\,\, i=1,2.\]
$\frac{\log |X_1|}{\log K}$ is a Poisson random variable with parameter $\lambda$, and moving any of the four primary lattice directions is equally likely. $S_t$ is strongly aperiodic recurrent and has infinite range, $E(|X_1|^m) < \infty$ for all $m < \infty$ (and, in particular, cov$(|X_1|) = \Gamma = \frac{1}{2} e^{(K^2 - 1)\lambda} I$), and every large jump causes a landing in a new copy of $A$. The only way to escape $\pi_K \inv A = \hp \inv \hat{A}$ is a step of size $K^0=1$.
\end{example}
\end{comment} % END COMMENTED-OUT EXAMPLE

Computational bounds on $\EV^{\hat x}(\hTDnc)$, by \eqref{eqn:Escape} and \eqref{eqn:EscapeDiscExp}, are 
\begin{equation} \label{eqn:EscapeDiscExpTorus}
\frac{n^2 - |x|^2}{\gamma^2} \leq \EV^{\hat x}(\hTDnc) \leq \frac{n^2 - |x|^2}{\gamma^2} + 2n + 1 + O(K^{-M} n^4).
\end{equation}

\begin{comment} % BEGIN ANOTHER COMMENTED-OUT EXAMPLE
\begin{example} \label{eg:LazySRW}
Define the $\varepsilon$-lazy simple random walk on $\Z^d$, for $0 \leq \varepsilon < 1$, to be the walk with steps 
$p_1(e_j) = p_1(-e_j) = \frac{1 - \varepsilon}{2d}$, $j=1,...,d$; $p_1(0) = \varepsilon$,
\ie the walk stands still for a step with probability $\varepsilon$, and acts ``simply'' otherwise. Then $\Gamma = \left(\frac{1-\varepsilon}{d}\right)I$, and so for $d=2$, $\EV^{\hat x}(\hTDnc) = \frac{n^2 - |x|^2}{1-\varepsilon} + O(n)$.
\end{example}
\end{comment} % END ANOTHER COMMENTED-OUT EXAMPLE

We will next see that, from inside a disc, the probability of hitting the center before escaping is nearly the same on the torus as on the plane. Recall that, for $\hat{x} \in \ZK$, $x := \pi_K \hp \inv \hat x$.

\begin{lem} \label{lem:HitZeroFirst}
For all $\hat{x} \in \hDn$ and $n$ sufficiently large with $2n < K/2$,
\begin{equation} \label{eq:HitZeroFirst}
P^{\hat{x}}(T_{\hat{0}} < \hTDnc) = P^{x}(T_0 < T_{D(0,n)^c}) + O(K^{-M} n^2).
\end{equation}
\end{lem}

\pf The event $\{T_{\hat{0}} < \hTDnc\}$ can occur in two ways:
\begin{itemize}
\item The walk hits $\hat{0}$ after a small jump, never leaving the disc. This is equivalent to the planar event $\{T_{\hat{0}} = T_0 < T_{D(0,n)^c}\}$.
\item The planar walk (wlog starting from $\pi_K x$) does not hit $0$, and exits $D(0,n)$ via a targeted jump into $\hp \inv \hDn$. It may do this multiple times before finally hitting $\pi_K \inv 0$ (via a small or large jump).
\end{itemize}
We can represent this event as the disjoint union
\begin{align*}
\{T_{\hat{0}} < \hTDnc\} = \,\,\, & \{T_{\hat{0}} = T_0 < T_{D(0,n)^c} \leq \hTDnc\} \\
 & \sqcup \{T_{D(0,n)^c} < T_{\hat{0}} < \hTDnc\}. 
\end{align*}
The first case contains $\{T_0 < T_{D(0,n)^c}\}$, so a lower bound on the toral probability is the planar result. An upper bound on the second case is found in the event $\{T_{D(0,n)^c} < \hTDnc\}$, which by (\ref{eqn:ProbTorusExit}) is rare. Hence,
\[ \begin{array}{lll}
P^x(T_0 < T_{D(0,n)^c}) & \leq & P^{\hat x}(T_{\hat{0}} < \hTDnc)\\
 & \leq & P^x(T_0 < T_{D(0,n)^c}) + P^x(T_{D(0,n)^c} < \hTDnc)\\
 & \leq & P^x(T_0 < T_{D(0,n)^c}) + O(K^{-M} n^2). \qed
\end{array}\]

Finally, we calculate bounds for hitting time probabilities of a small disc around zero before escaping the $n$-disc. Let $\rho(\hat{x}) := n-|\hat{x}|$ be the distance between $\hat{x}$ and $\hp(D(0,n))$. 

\begin{lem}
Let $0 < \delta < \varepsilon < 1$. Then there exist $0 < c_1 < c_2 < \infty$ such that for all $\hat x \in \hDn \setminus \hp(D(0,\varepsilon n))$, for $n$ sufficiently large,
\begin{equation} \label{eq:InnerHitBounds}
c_1 \frac{\rho(\hat x) \vee 1}{n} \leq P^{\hat x}(T_{\hp(D(0,\delta n))} < \hTDnc) \leq c_2 \frac{\rho(\hat x) \vee 1}{n}.
\end{equation}
\end{lem}

\pf From (\ref{eq:HittingTimeComp1}), it is clear that $T_{D(0,n)^c} \leq \hTDnc$. Note that $T_{D(0,\delta n)} < T_{D(0,n)^c}$ only if the walk enters $D(0,\delta n)$ via a small jump (of distance no more than $(1-\delta)n$), and in this case $T_{D(0,\delta n)} = T_{\hp(D(0,\delta n))}$.
A large jump automatically causes planar exit of $D(0,n)$, regardless of where in the torus the walk lands. Breaking down the sets of paths involved, we have the planar case
\[\{T_{D(0,\delta n)} < T_{D(0,n)^c}\} = \{T_{\hp(D(0,\delta n))} = T_{D(0,\delta n)} < T_{D(0,n)^c} \leq \hTDnc\}\]
which covers all small-jump entrances to $\hp\inv\hp(D(0,\delta n))$; for the toral case, we have
\begin{align} 
 & \{T_{\hp(D(0,\delta n))} < \hTDnc\} \notag\\
 & = \,\,\, \{T_{\hp(D(0,\delta n))} = T_{D(0,\delta n)} < T_{D(0,n)^c} \leq \hTDnc\} \label{eq:AnnulusSets} \\
 &  \,\,\, \sqcup  \,\,\, \{T_{D(0,n)^c} \leq T_{\hp(D(0,\delta n))} < \hTDnc, T_{D(0,\delta n)}\}, \notag
\end{align}
where the second set contains all paths where a large jump occurs at or before entry to the inner disc. Hence,
\[P^x(T_{D(0,\delta n)} < T_{D(0,n)^c}) \leq P^{\hat x}(T_{\hp(D(0,\delta n))} < \hTDnc),\]
and so we get the lower bound from \cite[Lemma 2.1]{BRFreq}.

The upper bound simply bounds the second set in \eqref{eq:AnnulusSets}. By \eqref{eqn:ProbTorusExit}, 
\begin{align*}
P^{\hat x}(T_{\hp(D(0,\delta n))} < \hTDnc) & \leq P^x(T_{D(0,\delta n)} < T_{D(0,n)^c}) \\
 & \,\,\, + P^x(T_{D(0,n)^c} < \hTDnc) \\
 & \leq P^x(T_{D(0,\delta n)} < T_{D(0,n)^c}) + O(K^{-M} n^2),
\end{align*}
and the error term is absorbed by the upper bound on $P^x(T_{D(0,\delta n)} < T_{D(0,n)^c})$. $\qed$

\subsection{Internal Green's function}

Here we will examine \emph{internal} Green's functions on the plane (\ie from inside a disc; Green's functions external to a disc will be analyzed in Section \ref{ch:Entry}). 
We extend some results of \cite{LawSRW} for symmetric random walks on $\Z^2$ to projections of these random walks onto $\Z_K^2$.

We define the Green's function in the usual way for $\hat x, \hat y \in \hp(A) = \hat{A} \in \ZK$ to be, in comparison to \eqref{eq:GreenDefn},  
\begin{equation} \label{eqn:Green}
\hat{G}_{\hp(A)}(\hat x,\hat y) := \sum_{j=0}^{\infty} P^{\hat x}(\hat S_j = \hat y; j < T_{\hp(A^c_K)}) 
\end{equation}
and 0 else. In the planar case, the stopping time $T_{A^c}$ for a bounded set $A$ has a clear meaning, as a sufficiently large jump (one with magnitude greater than the diameter of $A$, for example) will certainly exit $A$. Jumps targeting $A$ land, in $\Z^2$, in $\pi_K\inv A = \hp\inv\hat{A}$; on $\ZK$, they land in $\hat{A}$. This means that planar estimates must be adjusted to reach similar results on the torally-projected walk. 

Please note that \eqref{eqn:Green} is different from the planar Green's function on the periodic planar set $\pi_K \inv A$: 
\begin{equation} \label{eqn:GreenPT}
G_{\pi_K \inv A}(x,y) := \sum_{j=0}^{\infty} P^x(S_j = y; j < T_{\pi_K \inv (A^c_K)}), \,\, x, y \in \pi_K \inv A.
\end{equation}
We will explore this distinction in Section \ref{ch:Entry}.

Note that $S_j \in \hp\inv\hat{S}_j$ for every $j$. By \eqref{eq:HittingTimeComp1} it is clear that planar escape happens at or before toral escape. Hence, the number of planar visits is less than or equal to the number of toral visits; for any $x,y \in A \subset \pi_K \Z^2$, %such that diam($A$) $< K$ (\ie fits inside one copy of $\Z^2_K$), 
\begin{align} 
G_A(x,y) & = \sum_{j=0}^{\infty} P^x(S_j = y; \, j < T_{A^c}) \notag\\
 & = \sum_{j=0}^{\infty} P^x(S_j \in \pi_K\inv y; \, j < T_{A^c}) 
 = \sum_{j=0}^{\infty} P^{\hat x}(\hat{S}_j = \hat{y}; \, j < T_{A^c}) \label{eqn:IntGreenIneq}\\
 & \leq \sum_{j=0}^{\infty} P^{\hat x}(\hat{S}_j = \hat{y}; \, j < T_{\hp(A^c_K)}) = \hat{G}_{\hp(A)}(\hat x, \hat y), \notag
\end{align}
where equality occurs between the first and second lines because, of all the copies of $y$ in $\pi_K\inv y$, only the primary copy $y = \pi_K y$ can be hit before the planar escape time $T_{A^c}$.

We start by giving bounds on the number of visits to $\hat{0}$ before escaping a disc. %\index{Green@$G_{D(0,n)}$}
\begin{lem} \label{lem:GreenZero}
For $n$ sufficiently large (with $2n < K/2$),
\begin{equation} \label{eqn:GreenZero}
\hGDn(\hat 0,\hat 0) = G_{D(0,n)}(0,0)[1 + O(K^{-M} n^2)].
\end{equation}
\end{lem}

\pf Our lower bound $G_{D(0,n)}(0,0) \leq \hGDn(\hat 0,\hat 0)$ is clear from \eqref{eqn:IntGreenIneq}. To achieve the upper bound, first decompose the count, noting that the toral event $\{\hat S_j = \hat 0; j < T_{\hDnc}; T_{D(0,n)^c} = T_{\hDnc}\}$ equals the planar event $\{S_j = 0; j < T_{\Dnc}; T_{D(0,n)^c} = T_{\hDnc}\}$. Applying the strong Markov property at $T_{D(0,n)^c}$, 
\begin{align}
\hGDn(\hat 0,\hat 0) = & \,\, \sum_{j=0}^{\infty} P^{\hat 0}\left(\hat S_j = \hat 0; j < T_{\hDnc}\right) \label{eq:hGDn00bound1}\\
 = & \,\, \sum_{j=0}^{\infty} P^{\hat 0}\left(\hat S_j = \hat 0; j < T_{\hDnc}; T_{D(0,n)^c} = T_{\hDnc}\right) \notag\\
& + \sum_{j=0}^{\infty} P^{\hat 0}\left(\hat S_j = \hat 0; j < T_{\hDnc}; T_{D(0,n)^c} < T_{\hDnc}\right) \notag\\
 = & \,\, G_{\Dn}(0,0) + \EV^0\left[\hGDn(\hat{S}_{T_{D(0,n)^c}},\hat 0); \,\, T_{D(0,n)^c} < \hTDnc\right], \notag
\end{align}
where, on $\{T_{D(0,n)^c} < \hTDnc\}$, $\hat{S}_{T_{D(0,n)^c}}$ is the point in $\hDn$ that our walk lands once escaping the planar disc $D(0,n)$ via a targeted jump into a copy.

By (\ref{eqn:ProbTorusExit}), we know $P^0(T_{D(0,n)^c} < \hTDnc) = O(K^{-M} n^2)$. %Since we are on the torus, we can bound $\hGDn(\hat{S}_{T_{D(0,n)^c}},\hat 0)$ by the maximum over all $\hat x \in \hDn$. 
The strong Markov property applied at $T_0$ gives us the planar equality
\begin{equation} \label{eq:GreenXZeroPlanar}
G_{D(0,n)}(x,0) = P^x(T_0 < T_{D(0,n)^c}) \, G_{D(0,n)}(0,0)
\end{equation}
which implies $G_{D(0,n)}(x,0) \leq G_{D(0,n)}(0,0)$ for all $x \in D(0,n)$. This equality has a clear analog on the torus, by applying the strong Markov property at $T_{\hat 0}$:
\begin{equation} \label{eq:GreenTorusXZero}
\hGDn(\hat x,\hat 0) = P^{\hat{x}}(T_{\hat{0}} < \hTDnc) \, \hGDn(\hat 0,\hat 0),
\end{equation}
which, with \eqref{eqn:EscapeDiscExpTorus} implies, for all $\hat x \in \hDn$,
\begin{equation} \label{eq:GreenTorusXZeroIneq}
\hGDn(\hat x,\hat 0) \leq \hGDn(\hat 0,\hat 0) \leq \EV^{\hat{0}}(\hTDnc) < \infty.
\end{equation}
Thus, $\max_{\hat x \in \hDn} \hGDn(\hat{S}_{T_{D(0,n)^c}},\hat 0) = \hGDn(\hat 0,\hat 0)$, and by combining \eqref{eq:hGDn00bound1}, \eqref{eq:GreenTorusXZeroIneq}, and \eqref{eqn:ProbTorusExit}, we have 
\[\hGDn(\hat 0,\hat 0) \leq G_{D(0,n)}(0,0) + O(K^{-M} n^2) \, \hGDn(\hat 0,\hat 0),\]
which, when substituted back into itself gives, for some $c < \infty$,
\[\begin{array}{ll}
\hGDn(\hat 0,\hat 0) & \leq G_{D(0,n)}(0,0) + \sum_{j=1}^{\infty} (cK^{-M} n^2)^j \, G_{D(0,n)}(0,0)\\
 & = G_{D(0,n)}(0,0)\left(1 + O(K^{-M} n^2)\right). \qed
\end{array}\]

Define the \emph{potential kernel} for $X$ on $\Z^2$ as follows: for $x \in \Z^2$, \index{axy@$a(x)$} 
\begin{align}
a(x) := \lim_{n \to \infty} \sum_{j=0}^n [p_j(0) - p_j(x)]. \label{eq:axyDefn}
\end{align}
Combining the generality of rotation of \cite[Ch. III, Sec. 12, P3]{Spitzer} and \cite[Theorem 4.4.6]{LawSRW} and the infinite-range argument of \cite[Prop. 9.2]{BRFreq} gives, for covariance matrix $\Gamma$ and norm $\mc{J}^*(x) := |x \cdot \Gamma^{-1} x|$, as $|x| \to \infty$, 
\begin{align}
a(x) = \frac{2}{\pi_{\Gamma}} \log\mc{J}^*(x) + C(p_1) + o(|x|^{-1}), \label{eq:axyGeneralCalc}
\end{align}

where $C(p_1)$ is a constant depending on $p_1$ but not $x$, and $\pi_{\Gamma} = 2\pi\sqrt{\det \Gamma}$. For $\Gamma = cI$, this reduces to 
\begin{align}
a(x) & = \frac{1}{c \pi} \log\left(\frac{|x|}{\sqrt{c}}\right) + C(p_1) + o(|x|^{-1}) \notag\\
 & = \frac{1}{c \pi} \log|x| + C'(p_1) + o(|x|^{-1}) , \label{eq:axyCI}
\end{align}
where $C'(p_1) = C(p_1) - \frac{1}{2c\pi}\log c$. For simple random walk on $\Z^2$, $c=\frac{1}{2}$, and so this is, from \cite[Theorem 4.4.4]{LawSRW}, 
\begin{align}
a(x) & = \frac{2}{\pi} \log|x| + \frac{2\gamma + \log 8}{\pi} + o(|x|^{-1}), \label{eq:axy}
\end{align}
where $\gamma$ is Euler's constant. From here on, we will write \eqref{eq:axyCI} with the form 
\begin{align}
a(x) = \frac{2}{\pi_{\Gamma}} \log|x| + C'(p_1) + o(|x|^{-1}). \label{eq:axyCIgen}
\end{align}
By the argument in \cite[(2.8)-(2.12)]{BRFreq} (which calculates the overshoot estimate of $O(n^{-1/4})$ mentioned in the note after \cite[Prop. 6.3.1]{LawSRW}), and using \eqref{eq:axyCIgen}, we get a computational result for \eqref{eqn:GreenZero} if $\Gamma = cI$: 
\begin{align}
G_{D(0,n)}(0,0) & = \frac{2}{\pi_{\Gamma}}\log n + C' + O(n^{-1/4}) \label{eq:FP(2.13)} 
\end{align}
which implies the toral Green's function 
\begin{align}
\implies \hGDn(\hat 0,\hat 0) & = G_{D(0,n)}(0,0)(1 + O(K^{-M} n^2)) \notag\\
 & = \left(\frac{2}{\pi_{\Gamma}}\log n + C' + O(n^{-1/4})\right)(1 + O(K^{-M} n^2)) \notag\\
 & = \frac{2}{\pi_{\Gamma}}\log n + C' + O(n^{-1/4}). \label{eq:GreenZeroTorusVal}
\end{align}
For $x, y \in \Z^2$ such that $|x| \ll |y|$, we have, by a Taylor expansion around $y$, 
\begin{align} 
\log|y-x| & = \log|y| + O\left( \frac{|x|}{|y|} \right). \label{eq:logEstimate}
\end{align}
In particular, if $x \in D(0,2r)$ and $y \in D(0,R/2)^c$, with $R = 4mr$, we have 
\begin{align} 
\log|y-x| & = \log|y| + O\left( m^{-1} \right). \label{eq:logEstimate2}
\end{align}
Note that \eqref{eq:logEstimate} and \eqref{eq:logEstimate2} hold in the toral case without adjustment.

Let $\eta = \inf\{ t \geq 1: S_t \in \{0\} \cup D(0,n)^c\}$. Then, following the argument of \cite[(2.14)-(2.15)]{BRFreq}, since $a(x)$ is harmonic with respect to $p$, $a(S_{t \land \eta})$ is a bounded martingale. Hence, $|a(S_{t \land \eta})|^2$ is a submartingale, so $\EV|a(S_{t \land \eta})|^2 \leq \EV|a(S_{\eta})|^2 < \infty$, meaning $\{a(S_{t \land \eta})\}$ are uniformly integrable. Hence, by the optional stopping and bounded convergence theorems, \eqref{eq:axyCIgen}, and \eqref{eq:logEstimate2}, 
\begin{align*}
a(x) & = \lim_{t \to \infty} \EV^x(a(S_{t \land \eta})) = \EV^x(a(S_{\eta})) = \EV^x(a(S_{\eta}); \, S_{\eta} \neq 0) \\
 & = \sum_{y \in \bd D(0,n)_{n^{3/4}}} a(y) P^x(S_{\eta} = y) + \sum_{y \in D(0,n+n^{3/4})^c} a(y) P^x(S_{\eta} = y) \\
 & = \left( \frac{2}{\pi_{\Gamma}} \log n + C'(p_1) + o(|x|^{-1}) + O(n^{-1/4})\right) P^x(S_{\eta} \neq 0) + O(n^{-1/4}), 
\end{align*}
which, combining the error terms into $O(|x|^{-1/4})$, matches \cite[Prop. 6.4.3]{LawSRW}: 
\begin{align}
P^x&(T_0 < T_{D(0,n)^c}) = P^x(S_{\eta} = 0) = 1 - \frac{a(x) - O(n^{-1/4})}{\frac{2}{\pi_{\Gamma}} \log n + C' + O(|x|)^{-1/4}} \label{eq:FP(2.16)}\\
 & = 1 - \frac{\frac{2}{\pi_{\Gamma}} \log |x| + C' + O(|x|^{-1/4})}{\frac{2}{\pi_{\Gamma}} \log n + C' + O(n^{-1/4})} 
 = \left( \frac{\log (n/|x|) + O(|x|^{-1/4})}{\log n}\right)(1 + O((\log n)^{-1})). \notag
\end{align}
With \eqref{eq:HitZeroFirst}, %\eqref{eq:FP(2.13)}, and \eqref{eq:GreenZeroTorusVal}, 
we move this to the torus: 
\begin{align}
P^{\hat x}(T_{\hat{0}} < \hTDnc) & = \frac{\log(n/|\hat{x}|) + O(|\hat{x}|^{-1/4})}{\log(n)}\bigg(1 + O((\log n)^{-1})\bigg)   + O(K^{-M} n^2) \notag\\
 & = \frac{\log(n/|\hat{x}|) + O(|\hat{x}|^{-1/4})}{\log(n)}\bigg(1 + O((\log n)^{-1})\bigg). \label{eq:ProbZeroBeforeDisc}
\end{align}

Next, we examine $\hat{x} \in \hDR \setminus \hDr$. By the fact that a large targeted jump may land a planar walk into $\hp\inv \hDr \setminus D(0,r)$ (the set of any copy of $D(0,r)$ that is not the primary copy), we may transfer the planar results \cite[(2.20), (2.21)]{BRFreq} 
\begin{align} 
P^x(T_{D(0,r)} > T_{D(0,R)^c}) & = \frac{\log(|x|/r) + O(r^{-1/4})}{\log(R/r)} \label{eq:FP(2.20)} \\
P^x(T_{D(0,r)} < T_{D(0,R)^c}) & = \frac{\log(R/|x|) + O(r^{-1/4})}{\log(R/r)} \label{eq:FP(2.21)}
\end{align}
uniformly for $r < |x| < R$ to the toral results 
\begin{align}
P^{\hat x}(\hTDr > \hTDRc) & = \frac{\log(|\hat{x}|/r) + O(r^{-1/4})}{\log(R/r)} + O(K^{-M} R^2) \notag\\
 & = \frac{\log(|\hat{x}|/r) + O(r^{-1/4})}{\log(R/r)} \label{eq:ToralGamblersSuccess}  \\
P^{\hat x}(\hTDr < \hTDRc) & = \frac{\log(R/|\hat{x}|) + O(r^{-1/4})}{\log(R/r)} + O(K^{-M} R^2) \notag\\
 & = \frac{\log(R/|\hat{x}|) + O(r^{-1/4})}{\log(R/r)}. \label{eq:ToralGamblersRuin}
\end{align}

By \eqref{eq:GreenXZeroPlanar}, \eqref{eq:FP(2.13)}, \eqref{eq:GreenZeroTorusVal}, \eqref{eq:GreenTorusXZero}, \eqref{eq:FP(2.16)}, and \eqref{eq:ProbZeroBeforeDisc}, we get as corollaries calculations and bounds for $G_{D(0,n)}(x,0)$,  $\hGDn(\hat x,\hat 0)$, $G_{D(0,n)}(x,z)$, and $\hGDn(\hat x,\hat z)$: for $x \in \Dn$ and $\hat x \in \hDn$, %, matching \cite[(2.18)]{BRFreq}: \index{Green@$\hGDn$}
for some $C = C(p) < \infty$, 
\begin{align} 
G_{D(0,n)}(x,0) & = P^x(T_0 < T_{\Dnc}) \, G_{\Dn}(0,0) \notag\\
 & = \frac{\log(n/|x|) + O(|x|^{-1/4})}{\log(n)} \left( 1 + O((\log n)^{-1} \right) \left( \frac{2}{\pi_{\Gamma}} \log n + C' + O(n^{-1/4})\right) \notag\\ 
 & = \frac{2}{\pi_{\Gamma}}\log \bigg(\frac{n}{|x|}\bigg) + C + O(|x|^{-1/4}), \label{eq:FP(2.18)} 
\end{align}
\begin{align}
\hGDn(\hat x,\hat 0) & = \frac{2}{\pi_{\Gamma}}\log \bigg(\frac{n}{|\hat x|}\bigg) + C + O(|\hat x|^{-1/4}) \label{eq:GreenXZeroTorusCalc} \\
G_{D(0,n)}(x,z) & \leq G_{D(x,2n)}(0,z-x) \leq c \log n. \label{eq:FP(2.19)} \\
\hGDn(\hat x,\hat z) & = G_{D(0,n)}(x,z) + O(K^{-M} n^2 \log n) \leq c \log n.\label{eqn:GreenXZBounds}
\end{align}

Finally, we have the following result paralleling \eqref{eq:InnerHitBounds}. Recall that $\rho(\hat{x}) = n - |\hat{x}|$.
\begin{lem} \label{lem:GreenFP2.2}
For any $0<\delta<\varepsilon<1$ we can find $0<c_{1}<c_{2}<\infty$, such that for  all $\hat x \in \hDn \setminus \hp(D(0,\varepsilon n))$, $\hat y \in \hp(D(0, \delta n))$ and all $n$ sufficiently large such that $2n < K/2$,
\begin{equation} \label{eq:GreenFP2.2G}
c_{ 1}\frac{\rho( \hat x)\lor 1}{n}\leq \hat{G}_{\hDn}(\hat y,\hat x)\leq c_{2}\frac{\rho( \hat x)\lor 1}{n}.
\end{equation}
\end{lem}

\pf We follow our standard technique. \cite[Lemma 2.2, (2.38)]{BRFreq} gives bounds for the planar random walk's Green's function. The  toral version of this Green's function has a lower bound of the planar version: 
\begin{equation} \label{eq:GreenFP2.2GPlanar}
c_{1}\frac{\rho( x)\lor 1}{n}\leq G_{D(0, n) }( y,x)\leq \hGDn(\hat y,\hat x).
\end{equation}
For the upper bound, use \eqref{eqn:GreenXZBounds} and \cite[Lemma 2.2]{BRFreq} again, for some constant $c_2$:
\[\hGDn(\hat y,\hat x) \leq G_{D(0, n)}(y,x) + O(K^{-M} n^2 \log n) \leq c_2 \frac{\rho(x) \lor 1}{n}. \qed\]

\subsection{Local Time}

The \emph{local time} of a point $x \in \Z^2$ up to time $t$ is the amount of time the walk spends at $x$ up to time $t$:
\begin{equation} \label{eqn:LocalTime}
L_t^x := \sum_{j=0}^t 1_{x}(S_j).
\end{equation}
Likewise, the \emph{toral local time} of a point $\hat x \in \Z^2_K$ is simply the toral version of the planar local time: 
\begin{equation} \label{eqn:ToralLocalTime}
L_t^{\hat x} := \sum_{j=0}^t 1_{\hat x}(\hp S_j).
\end{equation}
It should be clear that, since $\hat x$ represents an infinite set of copies of $x \in \Z^2$, $L_t^{\hat x} \geq L_t^{x}$ for every $t$. Also, the expected local time is equal to the Green's function up to time $t$ for the point $x$:
\begin{align}
\EV^y(L_t^x) & = \EV^y \left(\sum_{j=0}^t 1_{x}(S_j)\right) = \sum_{j=0}^t P^y(S_j=x) = G_t(y,x)  \label{eq:ExpLocalTimeGreen} \\
\EV^{\hat y}(L_t^{\hat x}) & = \EV^{\hat y} \left(\sum_{j=0}^t 1_{\hat x}(\hp S_j)\right) = \sum_{j=0}^t P^{\hat y}(\hp S_j=\hat x) = {\hat G}_t(\hat y, \hat x).  \label{eq:ExpLocalTimeGreenToral}
\end{align}
Building on \cite[Section 4]{BRFreq}, we establish the straightforward relationship between local times and Green's functions on the torus-projected walk, and give bounds on the tail probabilities of a local time.

\begin{lem} \label{lem:TorusExpLocalTime}
For $\hat x \in \Z^2_K$, $|\hat x| < n < K/2$, if $L^{\hat{0}}_{T_{\hDnc}}$ is the toral local time at zero before escaping the toral $n$-disc, then its expectation is
\begin{equation} \label{eq:TorusELTG}
\EV^{\hat x}\left(L^{\hat{0}}_{T_{\hDnc}}\right) = \hat{G}_{\hDn}(\hat x,\hat 0).
\end{equation}
Furthermore, for all $z \geq 1$,
\begin{equation} \label{eq:TorusPLTzG}
P^{\hat x}\left(L^{\hat{0}}_{T_{\hDnc}} \geq z\hat{G}_{\hDn}(\hat 0,\hat 0) \right) \leq c\sqrt{z} e^{-z},
\end{equation}
for some $c < \infty$ independent of $\hat x, z, n$.
\end{lem}

\pf It should be obvious from \eqref{eqn:ToralLocalTime} that \eqref{eq:TorusELTG} is just \eqref{eq:ExpLocalTimeGreen}. To show \eqref{eq:TorusPLTzG}, we need the following machinery: by the Strong Markov property, we have for any power $k$,
\begin{align*}
\EV^{\hat x}\left[\left(L^{\hat{0}}_{T_{\hDnc}}\right)^k\right] & = k! \EV^{\hat x}\left[\sum_{0 \leq j_1 \leq j_2 \leq
\ldots \leq j_k \leq T_{\hDnc}} \prod_{i=1}^k 1_{\{\hat{S}_{j_i}=0\}}\right]\\
 & = k! \EV^{\hat x}\left[\sum_{0 \leq j_1 \leq j_2 \leq \ldots \leq j_{k-1} \leq T_{\hDnc}} \prod_{i=1}^{k-1} 1_{\{\hat{S}_{j_i}=0\}} \hat{G}_{\hDnc}(\hat 0,\hat 0)\right]\\
 & = k \EV^{\hat x}\left[\left(L^{\hat{0}}_{T_{\hDnc}}\right)^{k-1}\right] \hat{G}_{\hDnc}(\hat 0,\hat 0).
\end{align*}

By induction on $k$,
\begin{equation} \label{eq:TorusELTkG}
\EV^{\hat x}\left[\left(L^{\hat{0}}_{T_{\hDnc}}\right)^k\right] = k!\hat{G}_{\hDn}(\hat x,\hat 0)\hat{G}_{\hDn}(\hat 0,\hat 0)^{k-1},
\end{equation}
matching \cite[(4.5)]{BRFreq}.

To prove (\ref{eq:TorusPLTzG}), use (\ref{eq:TorusELTkG}),
(\ref{eq:GreenTorusXZero}), and Chebyshev's inequality to obtain
\begin{equation} \label{eq:TorusPLTzGbdFactorial}
P^{\hat x}\left(L^{\hat{0}}_{T_{\hDnc}} \geq z\hat{G}_{\hDn}(\hat 0,\hat 0)\right) \leq \frac{k!}{z^k},
\end{equation}
then take $k=[z]$ and use Stirling's approximation. $\qed$

%END Ch2

%\include{Ch3} % results related to entrance times
%\section{From the outside, looking in} % Disc entry functions 
\section{Disc Entry} \label{ch:Entry}

Here we will examine paths starting outside a disc. Since, on $\Z^2$, 
\begin{equation} \label{eq:AnnDiscSubs2}
\bd D(0,n)_s \subset \left\{\begin{array}{c}
\pi_K \inv \bd D(0,n)_s \\
D(0,n+s)
\end{array}\right\} \subset \pi_K \inv D(0,n+s),
%\,\,\, \bd D(0,n)_s \subset D(0,n) \subset \hDn,
\end{equation}
then starting at any $y \in \pi_K \inv (D(0,n+s)^c \cap \pi_K \Z^2)$ (as seen in Figure \ref{fig:discs_annuli}) yields the \emph{disc entrance time inequalities}
\begin{equation} \label{eq:HittingTimeComp2}
T_{\pi_K \inv D(0,n+s)} \leq \left\{\begin{array}{c}
T_{\pi_K \inv \bd D(0,n)_s} \\
T_{D(0,n+s)}
\end{array}\right\} \leq T_{\bd D(0,n)_s}.
\end{equation}

These relationships will be exploited in this and the next section. 

\subsection{External Green's function}

To supplement the internal Green's functions of Section \ref{ch:Escape} are external Green's functions: those counting the number of visits to a point outside of a set before entering that set. Wlog $x$ and $D(0,n)$ are in the primary copy. We will find bounds on three different external Green's functions: 
\begin{center}
\begin{tabular}{c|c|c|c|l}
Green's function & scope & starting at & counts visits to & before...\\
\hline
$G_{D(0,n)^c}(x,y)$ & planar & $x$ & $y$ & $T_{D(0,n)}$ \\
$G_{\pi_K \inv (D(0,n)^c_K)}(x,y)$ & planar & $x$ & $y$ & $T_{\pi_K \inv D(0,n)} = \hTDn$ \\
$\hGDnc(\hat{x},\hat{y})$ & toral & $\hat{x}$ & $\hat{y}$ & $T_{\pi_K \inv D(0,n)} = \hTDn$
\end{tabular}
\end{center}

While proofs for both escape and entry rely on the strong Markov property, there is a distinctly different flavor in the approaches used to prove entrance results due to the ``targeting effect'' induced by the periodicity of the toral projection. 

Note that, similar to \eqref{eq:GreenXZeroPlanar}, for any $x, y \in D(0,n)^c$, by the symmetry of $G_A$ and the strong Markov property at $T_x$, 
\begin{equation} \label{eq:ExternalPlanarEasyBound}
G_{D(0,n)^c}(x,y) = P^y(T_x < T_{D(0,n)}) \, G_{D(0,n)^c}(x,x),
\end{equation}
so, assuming $|x| < |y|$, we only need $G_{D(0,n)^c}(x,x)$ for an upper bound. Fix $j > 2$ and let 
\begin{align*}
U_0 & = 0, \\
V_i & = \min\{ t > U_i : |S_t| < n \, \mbox{ or } \, |S_t| > |x|^j\},\\
U_{i+1} & = \min\{ t > V_i : S_t = x\}.
\end{align*}
$U_i$ is the $i$th visit to $x$ after visiting $D(0,n)$ or $D(0,|x|^j)^c$ (noting that there can be multiple visits to $x$ before $V_i$, but none in the interval $V_i \leq t < U_{i+1}$). Hence, 
\begin{align*}
G_{D(0,n)^c}(x,x) & = \EV^x \left(\sum_{i=1}^{\infty} 1_{\{x\}}(S_t) \, 1_{\{t < T_{D(0,n)}\}}\right)\\
 & \leq \sum_{i=1}^{\infty} \EV^x \left( \sum_{U_i < t < V_i} 1_{\{x\}}(S_t) \, ; \, U_i < T_{D(0,n)}\right)\\
 & \leq \sum_{i=1}^{\infty} \EV^x \left( G_{D(0,|x|^j)}(S_{U_i}, x) \, ; \, U_i < T_{D(0,n)}\right)
\end{align*}
by the strong Markov property at $U_i$. Since $S_{U_i} = x$, and setting $a_i = a_i(x) := P^x(U_i < T_{D(0,n)})$, we have 
\begin{align*}
G_{D(0,n)^c}(x,x) & \leq \sum_{i=1}^{\infty} \EV^x \left( G_{D(0,|x|^j)}(S_{U_i}, x) \, ; \, U_i < T_{D(0,n)}\right) 
 = G_{D(0,|x|^j)}(x, x) \left(\sum_{i=1}^{\infty} a_i\right).
\end{align*}
To sum the $a_i$, note that by strong Markov at $U_i$ again, 
\begin{align*}
a_{i+1} & \leq \EV^x \left( P^{S_{U_i}}(T_{D(0,|x|^j)^c} < T_{D(0,n)}) \, ; \, U_i < T_{D(0,n)} \right)\\
 & = P^x(U_i < T_{D(0,n)}) \, P^x(T_{D(0,|x|^j)^c} < T_{D(0,n)})\\
\implies a_{i+1} & \leq a_i P^x(T_{D(0,|x|^j)^c} < T_{D(0,n)}),
\end{align*}
where, by \eqref{eq:FP(2.20)}, and for sufficiently large $n$, 
\begin{align}
P^x(T_{D(0,|x|^j)^c} < T_{D(0,n)}) & = \frac{\log(|x|/n) + O(n^{-1/4})}{\log(|x|^j/n)} \notag\\
 & \leq \frac{1 + O(n^{-1/4})}{j} \leq \frac{2}{j}. \label{eq:PlanarExponentialGamblersRuin}
\end{align}
Hence, $a_{i+1} \leq \frac{2}{j} a_i$, which implies $a_i \leq (2/j)^i$, and so by \eqref{eq:FP(2.13)}, 
\begin{align} 
G_{D(0,n)^c}(x,x) & \leq G_{D(0,|x|^j)}(x, x) \sum_{i=1}^{\infty} \left(\frac{2}{j}\right)^i \notag\\
 & \leq \frac{2/j}{1 - 2/j} G_{D(0,|x|^j + |x|)}(0, 0) \label{eq:ExternalGreenUpperBound}\\
 & \leq \frac{2}{j-2} G_{D(0,2|x|^j)}(0, 0) \leq \frac{2j}{j-2} G_{D(0,2|x|)}(0, 0). \notag
\end{align}

Moving the external Green's function to the torus is trickier: we must examine the conflict between counting visits to an infinite number of planar copies of $x$ in $\hp\inv\hat{x} = \pi_K \inv x$ versus avoiding an infinite number of copies of $D(0,n)$ in $\hp\inv\hDn = \pi_K \inv \Dn$. Also, the size of $j$ (via $|\hat x|^j$) is restricted relative to $K$. We use the same argument as before, with the adjustment of \eqref{eq:ToralGamblersRuin} applied to the argument of \eqref{eq:PlanarExponentialGamblersRuin}, yielding 
\begin{align*}
P^{\hat x}(T_{\hp(D(0,|\hat x|^j)^c_K)} < \hTDn) & = P^x(T_{D(0,|x|^j)^c} < T_{D(0,n)}) + O(K^{-M} n^2) \notag \\
 & \leq \frac{2}{j} + O(K^{-M} n^2).
\end{align*}
This gives the toral analog $\hat{a}_i \leq \left(\frac{2}{j} + O(K^{-M} n^2)\right)^i$, and the sum $\displaystyle{\sum_{i=1}^{\infty} \hat{a}_i}$ gives the slightly different bound of $\displaystyle{\sum_{i=1}^{\infty} \hat{a}_i \leq \frac{2}{j-2} + O(K^{-M} n^2)}$, yielding the toral upper bound %\index{Green@$\hGDnc$}
\begin{align} 
\hGDnc(\hat x,\hat x) & \leq \hat{G}_{\hp(D(0,|x|^j))}(\hat x,\hat x) \left(\sum_{i=1}^{\infty} \hat{a}_i\right) \notag\\
 & \leq \left( \frac{2}{j-2} + O(K^{-M} n^2) \right) j \hat{G}_{\hp(D(0,2|x|))}(0,0). \label{eq:ToralExternalGreenUpperBound}
\end{align}

Easy lower bounds for both the planar and toral cases are found by merely considering the visits to $x$ in the disc $D(x,|x|-n)$, whose boundary rests just outside $D(0,n)$. This bound, \eqref{eq:ExternalGreenUpperBound}, and \eqref{eq:ToralExternalGreenUpperBound} give, for $j > 2$, 
\begin{align} 
G_{D(0,|x|-n)}(0,0) & \leq G_{D(0,n)^c}(x,x) \leq \frac{2j}{j-2} G_{D(0,2|x|)}(0, 0) \label{eq:ExternalGreenBounds}\\
\hat{G}_{\hp(D(0,|x|-n))}(\hat 0,\hat 0) & \leq \hGDnc(\hat x,\hat x) \notag\\
 & \leq \left( \frac{2}{j-2} + O(K^{-M} n^2) \right) j \hat{G}_{\hp(D(0,2|x|))}(\hat 0,\hat 0). \label{eq:ToralExternalGreenBounds}
\end{align}

Thus, for any $x,y \in \pi_K(\Dnc_K)$ \st $|x| \leq |y|$, by (\ref{eqn:GreenXZBounds}) and (\ref{eq:GreenZeroTorusVal}), \eqref{eq:ExternalGreenBounds} and \eqref{eq:ToralExternalGreenBounds} become the computational bounds
\begin{eqnarray} 
G_{D(0,n)^c}(x,y) & \leq & \frac{2j}{j-2} \left[ \frac{2}{\pi_{\Gamma}} \log(2|x|) + C + O(|x|^{-1/4})\right] \leq c_j \log |x|, \label{eqn:ExtGreenIneq}\\
\hGDnc(\hat x,\hat y) & \leq & \frac{2j}{j-2} \left[ \frac{2}{\pi_{\Gamma}} \log(2|x|) + C + O(|x|^{-1/4})\right] \leq \hat{c}_j \log |\hat x|, \label{eqn:ExtGreenIneqToral}
\end{eqnarray}
where $c_j, \hat{c}_j$ depend on $j>2$, $c_j \geq \hat{c}_j$, and in the toral case, such that $|\hat x| < (\frac{K}{2})^{1/j}$ (there is no such restriction on the planar case). From here on, we consider $j=3$.

For $\left( \frac{K}{2} \right)^{1/3} < |x| < \frac{K}{2}$, first note that $\log|x| \asymp \log K$. 
By \eqref{eqn:ExtGreenIneq} and the fact that $\pi_K \inv (\Dnc \cap \pi_K \Z^2) = \hp\inv\hDnc \subset D(0,n)^c$, on $\Z^2$ we have \index{Green@$G_{\hp\inv\hDnc}$}
\begin{align}
G_{\hp\inv\hDnc}(x,x) \leq G_{D(0,n)^c}(x,x) \leq c \log |x| \label{eq:GhatDcVsGDc}
\end{align}
for any $x \in \pi_K(\Dnc_K)$.
%(primary) $x \in \pi_K \hp\inv\hDnc$. 
We can relate 
%$G_{\hp\inv\hDnc}(x,x)$
$G_{\Dnc}(x,x)$ to $\hGDnc(\hat x,\hat x)$ by the following inductive strategy: define 
\begin{align*}
T_0 & = 0 \\
T_1 & = \inf \left\{ t > 0 \, | \, S_t \in \hp\inv\hat{x} \setminus \{x\} \right\} \\
\vdots \\
T_{j+1} & = \inf \left\{ t > 0 \, | \, S_t \in \hp\inv\hat{x} \setminus \{ S_{T_j}\} \right\}, \, j=1,2,3,\ldots 
\end{align*}
as the hitting times of distinct copies of $x$. Let $U = \hTDn$. Then 
\begin{align*}
\hGDnc(\hat x,\hat x) = \sum_{j=0}^{\infty} 
  \EV^x \left(\sum_{ T_j \leq t < T_{j+1} \land U } 1_{\{S_t = S_{T_j} \}} \right).
\end{align*}
By the Markov property, translation invariance for different points of $\hp\inv \hat x$, and \eqref{eq:GhatDcVsGDc}, 
\begin{align*}
 \quad \EV^x \left( \sum_{T_j \leq t < T_{j+1} \land U} 1_{\{S_t = S_{T_j}\}} \right) 
 & = \EV^x \left( \sum_{T_j \leq t < T_{j+1} \land U} 1_{\{S_t = S_{T_j}\}} , T_j < U \right) \\
 & = P^x(T_j < U) \, \EV^x \left( \sum_{0 \leq t < T_{1} \land U} 1_{\{S_t = x\}} \right) \\
% & \leq P^x(T_j < U) \, G_{\hp\inv\hDnc}(x,x) \\
 & \leq P^x(T_j < U) \, G_{\Dnc}(x,x).
\end{align*}
Hence, 
\begin{align} 
\hGDnc(\hat x,\hat x) \leq G_{\Dnc}(x,x) \sum_{j=0}^{\infty} P^x(T_j < U). \label{eq:ToralGPlanarG}
\end{align}
Let $\rho = P^x(T_1 < U)$. By the Markov property, $P^x(T_j < U) \leq \rho^j$, so since $\rho < 1$, by \eqref{eq:GhatDcVsGDc} and \eqref{eq:ToralGPlanarG}, we have, for some $C < \infty$, 
\begin{align} 
\hGDnc(\hat x,\hat x) \leq \frac{C \log K}{1 - \rho}. \label{eq:ToralGBigBound}
\end{align}
Examining $\rho$, note that 
\begin{align*}
\{T_1 < U\} & = \left\{ T_{\pi_K \inv x \setminus \{x\}} < \hTDn \right\} \\
 & \subset \left\{ T_{D(0,K)^c} < \hTDn \right\} \subset \left\{ T_{D(0,K)^c} < T_{D(0,n)} \right\},
\end{align*}
which, by \eqref{eq:FP(2.21)} gives the bound 
\begin{align} \label{eq:RhoUpperBound} 
\rho & \leq 1 -  P^x \left(T_{D(0,n)}< T_{D(0,K)^c}  \right) \leq 1 - \frac{ \log(K/|x|) + O(n^{-1/4}) }{ \log(K/n) } \leq 1 - \frac{c}{\log K} 
\end{align}
for some $c < \infty$. Therefore, for $(\frac{K}{2})^{1/3} \leq |x| < \frac{K}{2}$, combining \eqref{eq:ToralGBigBound} and \eqref{eq:RhoUpperBound} we have the upper bound 
\begin{align} \label{eq:ToralExternalGUpperBound2}
\hGDnc(\hat x,\hat x) \leq C (\log K)^2.
\end{align}

\begin{figure}[!ht]
  \centering
    \includegraphics[width=4in]{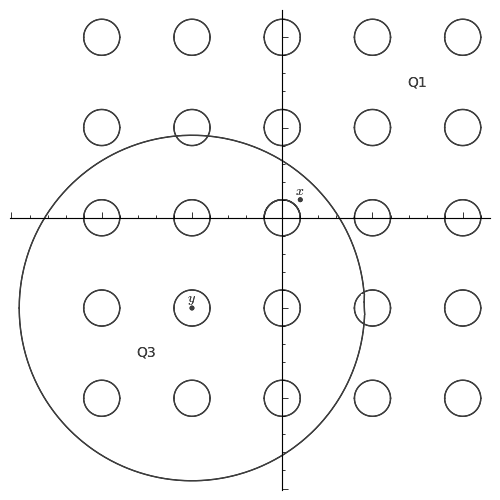}
  \parbox{4in}{
  \caption{$x \in Q1 \cap D(0,n)^c$, $y = (-K,-K)$. }
  \label{fig:external_green_planar}}
\end{figure}

Finally, for $|x| \geq \frac{K}{2}$, and considering $x = \pi_K x$, condition on the quadrant containing $S_{T_1}$ via a gambler's ruin from the opposite quadrant. For example, as in Figure \ref{fig:external_green_planar}, if $S_{T_1}$ is in Quadrant 1, let $y = (-K, -K)$ (opposite, in $Q3$). Then there exists $c$, $\left(\sqrt{2} + \frac{1}{2}\right)K \leq c < 2\sqrt{2}$ (say, $c=\sqrt{2} + \frac{1}{2}$ if $x \not \in Q1$, and $c=\frac{7}{4}\sqrt{2}$ if $x \in Q1$), such that, again by \eqref{eq:FP(2.21)}, we have the bound 
\begin{align} \label{eq:RhoConditionedUpperBoundFarAway}
P^x(T_1 < U \, | \, S_{T_1} \in Qi ) & \leq 1 - P^x \, \left( T_{D(y,n)} < T_{D(y,cK)^c} \right) \notag \\
 & \leq 1 - \frac{ \log(cK/|x-y|) + O(n^{-1/4}) }{ \log(cK/n) } \leq 1 - \frac{c'}{\log K} 
\end{align}
for some $c' < \infty$. 

%$\{T_1 < U \, | \, S_{T_1} \in Q1\}$

%All four quadrants share this estimate.
%By the periodicity of $\hat{x}$ and the symmetry of $\hat{p}_1$, the probability of $X_{T_1}$ being in any particular quadrant is roughly equal: $P^x(X_{T_1} \in Qi) = 1/4 + o(K^{-1})$, $i=1,2,3,4$. 
Thus, we have a bound on $\rho$ of 
\begin{equation} \label{eq:RhoUpperBoundFarAway}
\rho = \sum_{i=1}^4 P^x(T_1 < U \, | \, S_{T_1} \in Qi ) P(S_{T_1} \in Qi)\leq 1 - \frac{c^*}{\log K} 
\end{equation}
for some $c^* < \infty$. Therefore, by \eqref{eq:ToralGBigBound} and \eqref{eq:RhoUpperBoundFarAway}, \eqref{eq:ToralExternalGUpperBound2} holds for $|x| > \frac{K}{2}$. Collecting these results and applying them to the torus, we have proven 
\begin{lem} \label{lem:ExternalGreenToralComplete}
For $\hat x \in \hDnc$, 
\begin{equation} \label{eqn:ExternalGreenToralComplete}
\hGDnc(\hat x,\hat x) \leq \left\{ \begin{array}{ll}
C \log |\hat x| & n < |\hat x| < \left( \frac{K}{2} \right)^{1/3} \\
C\log^{2} |\hat x| & \left( \frac{K}{2} \right)^{1/3} \leq |\hat x|.
\end{array}\right.
\end{equation}
\end{lem}

\subsection{Disc entrance time}

We will now approach disc entrance times.
%, finding means and probability estimates for the inequalities in (\ref{eq:HittingTimeComp2}). 
Our first planar result mirrors (\ref{eqn:EscapeDiscExp}), with a very different end result, which is hinted by the first passage time result for SRW on $\Z$ (in, for example, \cite{ShreveVolI}).

\begin{lem} \label{lem:EnterDisc}
For any $y \in D(0,n)^c$, \index{expected@$\EV^y(T_{D(0,n)})$}
\begin{equation} \label{eqn:EnterDiscExp}
\EV^y(T_{D(0,n)}) = \infty.
\end{equation}
\end{lem}

\pf As in Lemma \ref{lem:EscapeDisc}, we use the martingale $M_t = |S_t|^2 - \gamma^2 t$, only this time we stop it at the stopping time $T_N := T_{D(0,n)} \land T_{D(0,N)^c}$, for $n \leq |y| < N$. Thus, the martingale stopped at $T_N \land k$ has expected value 
\[\EV^y(M_{T_N \land k}) = |y|^2 = \EV^y \left(|S_{T_N \land k}|^2 - \gamma^2 (T_N \land k)\right).\]
It is clear that $T_N \leq T_{D(0,N)^c}$, so by the argument given in the proof of Lemma \ref{lem:EscapeDisc}, we have $\EV^y \left( |S_{T_N \land k}|^2 \right) \to \EV^y \left( |S_{T_N}|^2 \right)$ and $\EV^y(T_N \land k) \to \EV^y(T_N)$ as $k \to \infty$. Hence, as in Lemma \ref{lem:EscapeDisc}, letting $\EV^y(M_{T_N \land k}) \to \EV^y(M_{T_N})$, we have 
\[\EV^y(T_N) = \frac{\EV^y \left(|S_{T_N}|^2\right) - |y|^2}{\gamma^2}.\]
Using \eqref{eq:FP(2.21)}, we decompose $\EV^y \left(|S_{T_N}|^2\right)$ and achieve the lower bound 
\begin{align*}
\EV^y \left(|S_{T_N}|^2\right) & = \EV^y \left(|S_{T_N}|^2 \, | \, T_N = T_{D(0,n)}\right) P^y(T_{D(0,n)} < T_{D(0,N)^c})\\
 & \,\,\, + \EV^y \left(|S_{T_N}|^2 \, | \, T_N = T_{D(0,N)^c}\right) P^y(T_{D(0,n)} > T_{D(0,N)^c})\\
 & \geq 0 + N^2 P^y(T_{D(0,n)} > T_{D(0,N)^c})\\
 & \geq \frac{N^2 \left[\log\left(|y|/n\right) + O(n^{-1/4})\right]}{\log\left(N/n\right)}
 \geq cN
\end{align*}
for some $c < \infty$. This gives us a lower bound on the expected entrance time of 
\[\EV^y(T_N) \geq \frac{c N - |y|^2}{\gamma^2},\]
which clearly goes to $\infty$ as $N \to \infty$. $\qed$

Next, we find finite bounds on the expected time to enter a toral disc. 
\begin{lem}
For any $n < \frac{K}{6}$ and $\hat y \in \hDnc$, there exists $c < \infty$ such that 
\begin{equation} 
\EV^{\hat y}(\hTDn) \leq \left\{ \begin{array}{ll}
c K^2 \log \left(n\right) & n < |\hat y| < n^2 \\
c K^2 \log \left(\frac{|\hat y|}{n}\right) & n^2 \leq |\hat y| < \left( \frac{K}{2} \right)^{1/3} \label{eq:ToralDiscEntryExpectedTime}\\
c K^2 (\log |\hat y|)^2 & \left( \frac{K}{2} \right)^{1/3} \leq |\hat y|.
\end{array}\right. 
\end{equation}
Also, we have the lower bound 
\begin{equation} 
\EV^{\hat y}(\hTDn) \geq \left\{ \begin{array}{ll}
\frac{(|\hat{y}|-n)^2}{\gamma^2} & |\hat{y}| < \frac{K}{3} \\
 \frac{c(K-n)^2}{\gamma^2} & |\hat y| \geq \frac{K}{3}
\end{array} \label{eq:ToralDiscEntryExpectedTimeLower} \right. 
\end{equation}
where $\gamma^2$ is as in the proof of Lemma \ref{lem:EnterDisc}.
\end{lem}

\pf For the upper bound, let $\hat y \in \hDnc$. We have the decomposition \index{expected@$\EV^{\hat y}(\hTDn)$}
\begin{equation*}
\EV^{\hat y}(\hTDn) = \sum_{\hat z \in \hDnc} \hGDnc(\hat y,\hat z),
\end{equation*}
which, for $|\hat y| > \left( \frac{K}{2} \right)^{1/3}$, by \eqref{eqn:ExternalGreenToralComplete} is clearly bounded above by $cK^2 (\log K)^2$. For closer $y$, we further decompose to 
\begin{align*}
\EV^{\hat y}(\hTDn) & = \sum_{\hat z \in \hp(D(0,|\hat y|)) \setminus \hDn} \hGDnc(\hat y,\hat z) \\
 & + \sum_{z \in \hp(D(0,|\hat y|)^c)} \hGDnc(\hat y,\hat z), 
\end{align*}
which, by \eqref{eqn:ExternalGreenToralComplete}, is bounded by 
\begin{align*}
\EV^{\hat y}(\hTDn) & \leq \sum_{\hat z: \, n \leq |\hat z| < |\hat y|} c \log |\hat z| \, + \, \sum_{\hat z: \, |\hat y| < |\hat z|} c \log |\hat y|\\
 & \leq c \int_0^{2\pi} \int_n^{|\hat y|} w \log w \, dw \, d\theta \, + \, c (K^2 - \pi |\hat y|^2) \log|\hat y|\\
 & \leq 2\pi c \left( \frac{1}{2} |\hat y|^2 \log|\hat y| - \frac{|\hat y|^2}{4} - n^2 \log n + \frac{n^2}{4}\right) \, + \, c(K^2 - \pi |\hat y|^2) \log|\hat y|\\
 & \leq c K^2 \log|\hat y| \leq c K^2 \log\left(\frac{|\hat y|}{n}\right) + cK^2 \log n - 2\pi c n^2 \log n.
\end{align*}
If $|\hat{y}| < \frac{K}{3}$, then 
$\EV^{\hat y}(\hTDn) \geq \EV^{\hat y}\left(T_{\hp(D(\hat{y}, |\hat{y}|-n)^c_K)}\right)$, 
and \eqref{eq:ToralDiscEntryExpectedTimeLower} follows directly from \eqref{eqn:EscapeDiscExpTorus}. The far-off $|\hat y| \geq \frac{K}{3}$ follows directly from the nearby case, since by the strong Markov property at $T_{\hp(D(0,K/3))}$, 
\begin{align*} 
\EV^{\hat y}(\hTDn) & \geq \EV^{\hat y}(T_{\hp(D(0,K/3))}; \, T_{\hp(D(0,K/3))} = T_{\hp(\bd D(0,K/6)_{K/6})})\\
 & \,\,\,  + \EV^{\hat y}( \EV^{\hat S_{T_{\hp(\bd D(0,K/6)_{K/6})}}}(\hTDn)) \geq \frac{c(K-n)^2}{\gamma^2}. \qed
\end{align*}

%\eqref{eq:ToralDiscEntryExpectedTime} hints at the cover time (and so, late points) results of Chapter \ref{ch:Late}. We will improve on these bounds in our discussion on excursions.

%END Ch3

%\include{Ch4} % results related to hitting annuli
%\section{Get in the ring} % Annulus functions 
\section{Annulus Entry} \label{ch:Annulus}

This section applies the results of \cite{CarThree} to the following partition of 
%\subsection{Application: Internal-External-Annulus} \label{sec:4.1}
%Let the following sets partition 
$\Z^2_K$ into a disc, an annulus around the disc, and the remainder of $\Z^2$, with $s \leq n < K \in \mathbb{N}$:
\begin{align*}
A = \hDn, \,\,\, B = \hDnsc, \,\,\, C = \hdDns.
\end{align*}
%\eqref{eq:BdEscapeEstimateToral} is an upper bound for \eqref{eq:psix} with our given $A$, $B$, $C$. 

\subsection{Hitting probabilities}

Starting from deep inside a disc, we first prove a bound on the probability of escaping the disc beyond an annulus outside it.
%state \cite[Lemma 2.3, (2.43)]{BRFreq}: 
%We improve on this bound for the toral version. 
\begin{lem} \label{lem:BdEscapeEstimateToral}
\begin{align} 
\sup_{x \in D(0,n/2)} & P^x(T_{\bd D(0,n)_s} > T_{D(0,n+s)^c}) & \leq c(s^{-M+2} \lor n^{-M+2}). \label{eq:BdEscapeEstimatePlanar} \\
\psi = \sup_{\hat x \in \hp(D(0,n/2))} & P^{\hat x}(\hTdDns > T_{\hDnsc}) & \leq c (s^{-M+2} \lor n^{-M+2}). \label{eq:BdEscapeEstimateToral}
\end{align}
\end{lem}

\pf To deal with targeting jumps that are possible on the torus, the proof of \eqref{eq:BdEscapeEstimateToral} is below. Replacing $\hat{p}_1(\hat y,\hat w)$ with $p_1(y,w)$, $\hat{G}_{\hDn}(\hat x,\hat y)$ with $G_{\Dn}(x,y)$,  \eqref{eqn:EscapeDiscExpTorus} with \eqref{eqn:EscapeDiscExp}, and \eqref{eq:GreenFP2.2G} with its planar version yields the proof for \eqref{eq:BdEscapeEstimatePlanar}.

For $\hat x \in \hp(D(0,n/2))$, we begin, by the last exit decomposition and \eqref{eqn:TargetedJumpEstimate}, with 
\begin{align*}
P^{\hat x}(\hTdDns & > T_{\hDnsc})
 = \sum_{\hat w \in \hDnsc} \sum_{\hat y \in \hDn}  \hat{G}_{\hDn}(\hat x,\hat y) \hat{p}_1(\hat y,\hat w) \notag\\
 & = \sum_{|\hat y| \leq 3n/4} \hat{G}_{\hDn}(\hat x,\hat y) \sum_{|\hat w| \geq n+s} ( p_1(  y,  w) + O(K^{-M})) \notag\\
 & \,\,\, + \sum_{3n/4 < |\hat y| < n} \hat{G}_{\hDn}(\hat x,\hat y) \sum_{|\hat w| \geq n+s} (  p_1(  y, w) + O(K^{-M})). \notag
\end{align*}
By \eqref{eq:GreenHit}, \eqref{eqn:EscapeDiscExpTorus}, and the facts that $|\hat y-\hat w| > n/4$ and $K > n$, the first sum has the bound 
\[c \sum_{|\hat y| \leq 3n/4} \hat{G}_{\hDn}(\hat x,\hat y) n^{-M} \leq c n^{-M+2}.\]
By \eqref{eq:GreenFP2.2G} and switching to a polar integral, the second sum is bounded by 
\begin{align*}
c \sum_{3n/4 < |\hat y| < n} \frac{n - |\hat y|}{n}(n+s-|\hat y|)^{-M}
 & \leq \frac{c}{n} \sum_{3n/4 < |\hat y| < n} \frac{n - |\hat y|}{(n+s-|\hat y|)^{M}} \\
 \leq \frac{c}{n} \sum_{3n/4 < |\hat y| < n} (n+s-|\hat y|)^{-M+1} 
 & \leq c \int_{3n/4}^{n} (n+s-u)^{-M+1} du \\
 & \leq c \int_0^{n/4} (s+v)^{-M+1} dv \leq cs^{-M+2}. \qed
\end{align*}

Note that for $\hat x \in \hDn$, by \eqref{eq:HittingTimeComp1},
\[ \{\hTdDns > \hTDnc\}^c = \{\hTdDns = \hTDnc\}.\]
Hence, provided $\hat x \in \hp(D(0,n/2))$, and $s \leq n$,  \eqref{eq:BdEscapeEstimateToral} is a bound for $\psi_{\hat x}$ from \cite[(7)]{CarThree}. { % ZZZ PREVIOUSLY \eqref{eq:psix}.
Also, \eqref{eq:ProbZeroBeforeDisc} and \eqref{eq:BdEscapeEstimateToral} gives us the chance of escaping a disc, into its $s$-annulus, before visiting its center:
\begin{align} \label{eq:FP(2.48)}
& P^{\hat x}(T_{\hat{0}} > \hTDnc; \, \hTDnc = \hTdDns) \notag\\
 & = 1 - \frac{ \log(n/|\hat x|) + O(|\hat x|^{-1/4}) }{ \log n }(1 + O((\log n)^{-1})  + O(s^{-M+2}).
\end{align}

By \eqref{eq:BdEscapeEstimateToral} and \cite[(9)-(10)]{CarThree}, 
for ${\hat x} \in \hp(D(0,n/2))$ and ${\hat y} \in \hDnsc$, 
\begin{align}
\rho_{\hat x} & = P^{\hat x}(T_{\hDnsc}, T_{\hDn}^* < T_{\hdDns}) \notag\\
& \leq c(s^{-M+2} \lor n^{-M+2}); \label{eq:rhoxBound} \\
\phi_{\hat y} & = P^{\hat y}(T_{\hDn}, T_{\hDnsc}^* < T_{\hdDns}) \notag\\
& \leq c(s^{-M+2} \lor n^{-M+2}). \label{eq:phixBound} 
\end{align}

Next, we find a %toral equivalent of \cite[Lemma 2.6]{BRFreq}, which 
bound for $\sigma_{\hat x}$ from \cite[(8)]{CarThree}.
\begin{lem} \label{lem:FPLemma2.6}
For $n$ sufficiently large, 
\begin{align}
\sigma & = \sup_{\hat x \in \hDnsc} P^{\hat x}(T_{\hDn} < T_{\hdDns}) \notag \\
 & \leq cn^2 \log(n)^2 (s^{-M} + n^{-M}). \label{eq:FP(2.71)}
\end{align}
\end{lem}
\pf Apply \eqref{eqn:ExternalGreenToralComplete} to the last exit decomposition to get 
\begin{align}
 & \sup_{\hat x \in \hDnsc} P^{\hat x}(T_{\hDn} < T_{\hdDns}) \label{eq:FP(2.72)}\\
\,\, = \,\,  & \sup_{\hat x \in \hDnsc} \sum_{\stackrel{\hat y \in \hp(D(0,n+s)^c)}{\hat w \in \hDn}} \hat{G}_{\hDnc}(\hat x,\hat y) \hat{p}_1(\hat y,\hat w) \notag\\
\,\,  \leq \,\, & c \log(n)^2 \sum_{n+s \leq |\hat y| < 2n} (|\hat y|-n)^{-M} + c \sum_{2n \leq |\hat y|} \log(|\hat y|)^2 (|\hat y|-n)^{-M} \notag\\
\,\,  \leq \,\, & cn^2 \log(n)^2 s^{-M} + c \log(n)^2 n^{-M+2}. \qed \notag
\end{align}

In particular, if $s=O(n)$, since $M=4+2\beta$, \eqref{eq:FP(2.71)} is bounded above by $cn^{-2}$, and if $s=O(\sqrt{n})$, \eqref{eq:FP(2.71)} is bounded above by $cn^{-\beta}$.

Combining \eqref{eq:ToralGamblersSuccess}  and \eqref{eq:BdEscapeEstimateToral}, we find the probability that, starting far from a small disc $\hDr$, the walk escapes a larger disc $\hDR$ before entering $\hDr$. If $r < R$ and $\hat x \in \hp(D(0,R/2))$,
% and $R^{-M+2} \log R < s^{-M+2}$, 
we have 
\begin{align} \label{eq:FP2.49analog}
& P^{\hat x}(\hTDRc < \hTDr; \hTDRc = T_{\hp(\bd D(0,R)_s)})\notag\\
& = \frac{\log(|\hat x|/r) + O(r^{-1/4})}{\log(R/r)} + O(s^{-M+2}).
\end{align}

To enter a disc, we first quote the planar result \cite[Lemma 2.4]{BRFreq}: if $s < r < R$ sufficiently large with $R \leq r^2$ we can find $c < \infty$ and $\delta > 0$ such that for any $ r < |x| < R$, 
\begin{align} \label{eq:FP(2.50)}
P^x(T_{D(0,r)} < T_{D(0,R)^c}; T_{D(0,r)} = T_{D(0,r-s)}) \leq cr^{-\delta} + cs^{-M+2}.
\end{align}
We see the same result on $\Z^2_K$, with an extra toral term (which is absorbed). 
\begin{lem} \label{lem:FPLemma2.4}
For the conditions listed above, 
\begin{align} 
P^{\hat x}(\hTDr & < \hTDRc; \hTDr = T_{\hp(D(0,r-s))}) \notag \\
 & \leq cr^{-\delta} + cs^{-M+2}. \label{eq:FP2.50analog}
\end{align}
\end{lem}

\pf Wlog we assume $\hat x \in \pi_K(D(0,R) \setminus D(0,r))$. Let 
\[\hat A := \{\hTDr < \hTDRc; \hTDr = T_{\hp(D(0,r-s))}\}.\]
Decompose $\hat A$ along the event $\{\hTDRc = T_{D(0,R)^c}\}$: by \eqref{eq:HittingTimeComp1} and Figure \ref{fig:discs_annuli}, 
\begin{align} \label{eq:EnterFarTorally}
P^{\hat x}(\hat A) & = P^{\hat x}(\hat A; \hTDRc = T_{D(0,R)^c}) + P^{\hat x}(\hat A; \hTDRc > T_{D(0,R)^c}). 
\end{align}
The first probability in \eqref{eq:EnterFarTorally} accounts for all walks that have no large jumps before the planar time $T_{D(0,r)}$, since 
\[\hat A \cap \{\hTDRc = T_{D(0,R)^c}\} \subset \{ \hTDr = T_{D(0,r)} \} \cap \{ \hTDr = T_{\hp(D(0,r-s))} \}.\]
Thus, \eqref{eq:FP(2.50)} bounds the first probability. The second probability is bounded by \eqref{eqn:ProbTorusExit}, which, as $cK^{-M} n^{2}$, is absorbed by $cs^{-M+2}$ since $s \leq R < K$.  $\qed$

We use \eqref{eq:FP2.50analog} along with \eqref{eq:ToralGamblersRuin} to get the toral gambler's ruin-via-annulus estimate:  
%(analog to \cite{BRFreq}, (2.66))
\begin{align} \label{eq:FP2.66analog}
& P^{\hat x}(\hTDr < \hTDRc; \hTDr = T_{\hp(\bd D(0,r-s)_s)}) \notag\\
& = \frac{\log(R/|\hat x|) + O(r^{-\delta})}{\log(R/r)} + O(s^{-M+2}).  %+ O(K^{-M+2+\beta}).
\end{align}

\subsection{Green's functions}

We start calculating bounds for the external Green's function with $\hat{x} \in \hp(D(0,n/2))$, $\hat{y} \in \hDn$: by 
\cite[(11)]{CarThree}, %\eqref{eq:GreenABUBa} with $A=\hDn$, 
\eqref{eqn:GreenXZBounds}, and \eqref{eq:rhoxBound}, 
\begin{align}
\hat{G}_{\hdDnsc}(\hat x,\hat y) & \leq \hat{G}_{\hDn}(\hat x, \hat y) + \frac{\rho_{\hat x}}{1 - \rho_{\hat y}} \hat{G}_{\hDn}(\hat y, \hat y). \label{eq:GhDnBound}  % \leq c \log n
\end{align}
In particular, if $\hat{y} = \hat{0}$ and $s = O(n)$, then $\rho_{\hat x} \leq c n^{-2}$ and by \eqref{eq:FP(2.18)}, 
\begin{align}
\hat{G}_{\hdDnsc}(\hat x,\hat 0) & \leq \hat{G}_{\hDn}(\hat x, \hat 0) + \frac{\rho_{\hat x}}{1 - \rho_{\hat 0}} \hat{G}_{\hDn}(\hat 0, \hat 0) \notag\\
\implies \hat{G}_{\hdDnsc}(\hat x,\hat 0) & = \frac{2}{\pi_{\Gamma}}\log\left(\frac{n}{|\hat x|}\right)  + C(\hat p_1) + O(|\hat x|^{-1/4}). \label{eq:GhDnBoundZero} 
\end{align}
By \cite[(12)]{CarThree}, %\eqref{eq:GreenABUBb}, 
\eqref{eqn:ExternalGreenToralComplete}, and \eqref{eq:phixBound}, for $\hat{x}, \hat{y} \in \hDnsc$, 
\begin{align}
\hat{G}_{\hdDnsc}(\hat x,\hat y) & \leq \hat{G}_{\hDnsc}(\hat x, \hat y)  + \frac{\phi_{\hat x}}{1 - \phi_{\hat y}} \hat{G}_{\hDnsc}(\hat y, \hat y)  \notag\\
 & \leq c (\log(|\hat{x}| \land |\hat{y}|))^2. \label{eq:GhDnscBound} 
\end{align}
Finally, for $\hat{x} \in \hp(D(0,n/2))$ and $\hat{y} \in \hDnsc$, by By \cite[(13)]{CarThree}, %\eqref{eq:GreenABUBab}, \eqref{eq:BdEscapeEstimateToral}, \eqref{eq:FP(2.71)}, and the above, 
\begin{align}
\hat{G}_{\hdDnsc}&(\hat x,\hat y) \label{eq:GreenbdDUB}\\
\leq & \min\left\{ \frac{\sigma_x}{1 - \rho_y} \hat{G}_{\hDn}(\hat x, \hat x), 
 \frac{\psi_x}{1-\phi_y} \hat{G}_{\hDnsc}(\hat y, \hat y) \right\} \notag\\
 \leq & \, c \, \min\left\{ n^2 (\log n)^3 (s^{-M} + n^{-M}), (\log(|\hat{y}|))^2  (s^{-M+2} \vee n^{-M+2}) \right\}. \notag
\end{align}
In particular, if $s=O(n)$, then in this case $\hat{G}_{\hdDnsc}(\hat x,\hat y) \leq c n^{-2}$, and if $s=O(\sqrt{n})$, the bound is $cn^{-\beta}$.

\subsection{Hitting times}

By (\ref{eq:HittingTimeComp2}) and (\ref{eqn:EnterDiscExp}), for $y \in D(0,n+s)^c \subset \Z^2$, the external planar annulus hitting time $E^y(T_{\bd D(0,n)_s}) = \infty$. Since, starting from inside the disc $x \in D(0,n)$, there is positive probability of hopping over an $s$-width annulus, then by the strong Markov property on $T_{D(0,n+s)^c}$, the internal planar annulus hitting time $E^x(T_{\bd D(0,n)_s}) = \infty$ as well. This is not the case for the toral analogues of these times. 

Torally, our walk can make small or targeted jumps before the disc escape time. To bound the annulus hitting times, we  employ \eqref{eqn:EscapeDiscExpTorus}, \eqref{eq:ToralDiscEntryExpectedTime}, and \cite[(21)]{CarThree}. %\eqref{eq:fAfB}. 
These yield, for some $c,c' < \infty$, 
\begin{align}
f_{\hDn}     = & \sup_{\hat x \in \hDn} E^{\hat x}(T_{\hDnc}) \leq cn^2, \label{eq:fD0n}\\
f_{\hDnsc} = & \sup_{\hat y \in \hDnsc} E^{\hat y}(T_{\hDns}) \leq c' (K \log K)^2. \label{eq:fD0nsc}
\end{align}
By  \cite[(22)-(23)]{CarThree}, %\eqref{eq:ExTCUBa}, \eqref{eq:ExTCUBb}, 
\eqref{eq:fD0n}, \eqref{eq:fD0nsc}, \eqref{eq:BdEscapeEstimateToral}, and \eqref{eq:FP(2.71)}, the expected annulus hitting time is bounded above: if $\hat{x} \in \hp(D(0,n/2))$ and $\hat{y} \in \hDnsc$, 
\begin{align} 
\EV^{\hat x}(T_{\hdDns}) & \leq \EV^{\hat x}(T_{\hDnc}) + \psi_{\hat x} \left[ \frac{f_{\hDnsc} + \sigma f_{\hDn}}{1 - \psi\sigma} \right] \notag\\
% & \leq \EV^{\hat x}(T_{\hDnc}) + c\left(s^{-M+2} \vee n^{-M+2}\right)\left((K \log K)^2 + (n \log n)^2 (s^{-M} + n^{-M})\right) \notag\\
 & \leq \EV^{\hat x}(T_{\hDnc}) + c\left(s^{-M+2} \vee n^{-M+2}\right)(K \log K)^2; \label{eq:InTbdDnsBounds}\\
\EV^{\hat y}(T_{\hdDns}) & \leq \EV^{\hat y}(T_{\hDns}) + \sigma_{\hat y} \left[ \frac{f_{\hDn} + \psi f_{\hDnsc}}{1 - \psi\sigma} \right] \notag\\
 & \leq c(K \log K)^2. \label{eq:ExTbdDnsBounds}
\end{align}
In particular, if $s, n=O(K)$, then for $K$ sufficiently large, note that by \eqref{eqn:EscapeDiscExpTorus}, %\eqref{eq:ToralDiscEntryExpectedTime}, and \eqref{eq:ToralDiscEntryExpectedTimeLower}, 
\begin{align*}
\EV^{\hat x}(T_{\hDnc}) = \frac{K^2 - |\hat x|^2}{\gamma^2} + O(K), %\,\,\,\, O(K^2) \leq \EV^{\hat y}(T_{\hDns}) \leq O( (K \log K)^2 ),
\end{align*}
%{\bf  the lower bound has conditions on $y$}
which, with $M = 4 + 2\beta$, reduces \eqref{eq:InTbdDnsBounds} to %and \eqref{eq:ExTbdDnsBounds}, yield the bounds 
\begin{align} 
\EV^{\hat x}(T_{\hdDns}) & = \left(1 + O(K^{-2-\beta})\right) \EV^{\hat x}(T_{\hDnc}). \label{eq:InTbdDnsBounds2}
% & = \frac{K^2 - |\hat x|^2}{\gamma^2} + O(K).  %\notag\\
%\EV^{\hat y}(T_{\hdDns}) & \leq \left(1 + O(K^{-\beta}(\log K)^2)\right) \EV^{\hat y}(T_{\hDns}). \label{eq:ExTbdDnsBounds2}
\end{align}

\singlespacing
%\include{Bibliography}
%\cleardoublepage


\addcontentsline{toc}{section}{References}

\begin{thebibliography}{120}

%\bibitem{AldP}
%Aldous, D (1989). 
%\emph{Probability Approximations via the Poisson Clumping Heuristic}.
%Springer, Berlin.

\bibitem{AF}
Aldous, D. and Fill, J.
\emph{Reversible Markov Chains and Random Walks on Graphs}.
\texttt{http://www.stat.berkeley.edu/$\sim$aldous/RWG/book.html}.
Monograph in preparation.  Last verified 13 July 2010.

\bibitem{BRFreq}
Bass, R. and Rosen, J. (2007).
Frequent points for random walks in two dimensions.
\emph{Electronic Journal of Probability} {\bf 12} 1-46.

\bibitem{CarThree}
Carlisle, Michael (2012).
On the Escape of a Symmetric Random Walk From Two Pieces of a Tripartite Set.
\texttt{arXiv:1209.1761 [math.PR]}.

%\bibitem{BH}
%Brummelhuis, M. and Hilhorst, H. (1991).
%Covering of a finite lattice by a random walk.
%\emph{Physica A.} {\bf 176} 387-408.

%\bibitem{DemboLN}
%Dembo, A. and Funaki, T. (2005)
%Lectures on probability theory and statistics: Ecole d'Et\'e de Probabilit\'es de Saint-Flour XXXIII - 2003.
%\emph{Lecture Notes in Mathematics} {\bf 1869}. Springer-Verlag Berlin Heidelberg.

%\bibitem{DemboICM}
%Dembo, A. (2006)
%Simple random covering, disconnection, late and favorite points.
%\emph{Proceedings of the International Congress of Mathematicians Madrid} {\bf 3} 535-558.

\bibitem{DPRZ2001}
Dembo, A., Peres, Y., Rosen, J., and Zeitouni, O. (2001).
Thick points for planar Brownian motion and random walks in two dimensions.
\emph{Acta Math.} {\bf 186} 239-270.

\bibitem{DPRZ2004}
Dembo, A., Peres, Y., Rosen, J., and Zeitouni, O. (2004).
Cover times for Brownian motion and random walks in two dimensions.
\emph{Annals of Mathematics} {\bf 160} 433-464.

\bibitem{DPRZ2006}
Dembo, A., Peres, Y., Rosen, J., and Zeitouni, O. (2006).
Late points for random walks in two dimensions. 
\emph{Annals of Probability} {\bf 34} 219-263.

%\bibitem{DiaconisGRPS}
%Diaconis, P. (1988).
%\emph{Group Representations in Probability and Statistics}.
%Inst. Math. Stat., Hayward.

%\bibitem{DiaconisCutoff}
%Diaconis, P. (1996).
%The cutoff phenomenon in finite Markov chains.
%In \emph{Group Representations in Probability and Statistics}.
%Inst. Math. Stat., Hayward.

\bibitem{Durrett}
Durrett, R. (2005).
\emph{Probability: Theory and Examples.} 3rd Edition.
Thomson - Brooks/Cole.

\bibitem{ET}
Erd\"os, P. and Taylor, S. J. (1960). 
Some problems concerning the structure of random walk paths.
\emph{Acta Math. Acad. Sci. Hungar.} 11, 137-162.

%\bibitem{FitzPit}
%Fitzsimmons, P. and Pitman, J. (1999).
%\emph{Kac's moment formula for additive functionals of a Markov process.}
%\emph{Stochastic Process. Appl.} {\bf 79} 117-134.

\bibitem{LawInt}
Lawler, G. (1991).
\emph{Intersections of Random Walks.}
Birkh\"auser, Boston.

\bibitem{LawCov}
Lawler, G. (1993).
On the covering time of a disc by a random walk in two dimensions.
In \emph{Seminar on Stochastic Processes 1992} 189-208. 
Birkh\"auser, Boston.

\bibitem{LawSRW}
Lawler, G. and Limic, V. (2010).
\emph{Random Walk: A Modern Introduction.}
%\texttt{http://www.math.uchicago.edu/~lawler/books.html}.
%Monograph in preparation.  Last verified 22 July 2009.
Cambridge studies in advanced mathematics, 123.
Cambridge University Press, New York.

\bibitem{LawPol}
Lawler, G. and Polaski, T. (1992).
Harnack inequalities and difference estimates for random walks with infinite range.
\emph{Journal of Theoretical Probability} {\bf 6} 781-802.

\bibitem{LPW}
Levin, D., Peres, Y., and Wilmer, E. (2008).
\emph{Markov Chains and Mixing Times}.
Amer. Math. Society.

%\bibitem{MR}
%Marcus, Michael B. and Jay Rosen. (2006).
%\emph{Markov Processes, Gaussian Processes, and Local Times}.
%Cambridge studies in advanced mathematics, 100.
%Cambridge University Press, New York.

\bibitem{RY}
Revuz, Daniel and Marc Yor. (2005).
\emph{Continuous Martingales and Brownian Motion}, 3rd Edition.
Springer Berlin Heidelberg New York.

\bibitem{RosenET}
Rosen, J. (2005).
A random walk proof of the Erd\"{o}s-Taylor Conjecture.
\emph{Periodica Mathematica Hungarica} {\bf 50} 223-245.

\bibitem{ShreveVolI}
Shreve, S. E. (2005).
\emph{Stochastic Calculus for Finance I: The Binomial Asset Pricing Model}.
Springer Science+Business Media, New York.

\bibitem{Spitzer}
Spitzer, F. (1976).
\emph{Principles of Random Walk}, Second Edition.
Springer, Princeton, NJ.

\bibitem{Wilf}
Wilf, H. (1989).
The editor's corner: the white screen problem.
In \emph{American Mathematical Monthly} {\bf 96} 704-707.

% book
%\bibitem{LABEL}
%Last, F. (YEAR).
%\emph{TITLE}.
%PUBLISHER LOCATION.
% I'd like to put ISBN as well.

% article
%\bibitem{LABEL}
%Last, F, and Last, F. (YEAR).
%ARTICLE TITLE.
%\emph{JOURNAL NAME} {\bf VOL\#} PA-GES.

% web page
%\bibitem{LABEL}
%Last, F.  
%``TITLE.''
%\texttt{URL}.
%Verified DATE.

\end{thebibliography}
\end{document}